\documentclass[a4paper, leqno]{article}
\usepackage{yhmath}
\usepackage{tgpagella}
\usepackage{euler}
\usepackage[mathscr]{euscript}
\usepackage[T1]{fontenc}
\usepackage{yfonts}
\usepackage[a4paper, total={6in, 8in}]{geometry}
\usepackage[utf8]{inputenc}
\usepackage{latexsym}
\usepackage{stmaryrd}
\usepackage{makeidx}
\usepackage{xr}
\usepackage{mathrsfs}
\usepackage{mathtools}
\usepackage{tikz-cd}
\usepackage{tikz}
\usepackage{pdflscape}
\usetikzlibrary{positioning}
\usepackage{float}  
\usepackage{amscd}
\usepackage{syntonly}
\usepackage{amssymb}
\usepackage{amsfonts}
\usepackage{amsmath}
\usepackage{amsthm}
\usepackage[english]{babel}
\usepackage{chngcntr}
\usepackage{xcolor}
\usepackage[numbers]{natbib}
\usepackage{adjustbox}
\usepackage{rotating}
\usepackage{quiver}
\usepackage{enumitem}
\usepackage{hyperref}
\usepackage{comment}
\usepackage{caption}
\usepackage{breqn}
\counterwithin{equation}{section}
\newtheorem{theorem}{Theorem}[section]

\newtheorem{definition}{Definition}[section]

\newtheorem*{note*}{Note}
\newtheorem*{lemma*}{Lemma}
\newtheorem*{examples*}{example}
\newtheorem*{example*}{example}
\newtheorem*{corollary*}{Corollary}
\newtheorem*{theorem*}{Theorem}
\newtheorem*{remark*}{Remark}
\newtheorem*{Reminder*}{Reminder}
\newtheorem*{Notation Convention*}{Notation Convention}
\newcommand{\dash}{\text{-}}
\newcommand{\la}[1]{\mathrm{#1}}
\newcommand{\ca}[1]{\mathcal{#1}}

\newcommand{\cat}[1]{\mathsf{#1}}
\newcommand{\go}[1]{\mathfrak{#1}}
\newcommand{\lu}[1]{\mathscr{#1}}

\def\dotminus{\mathbin{\ooalign{\hss\raise1ex\hbox{.}\hss\cr
  \mathsurround=0pt$-$}}}
  
\date{}
\author{Ali Hamad \footnote{Affiliation: University of Ottawa, Email: ahama099@uottawa.ca}}
\title{ Ultracategories as colax algebras for a pseudo-monad on $\mathsf{CAT}$}
\date{August 2026}

\begin{document}
\oddsidemargin 0in
\evensidemargin 0in
\topmargin=0in
\textwidth=6in
\textheight=8in
\maketitle
\begin{abstract}
 We show a result inspired by a conjecture by Shulman claiming that ultracategories as defined by Lurie   are normal colax algebras for a certain pseudo-monad on the category of categories $\cat{CAT}$. Such definition allows us to regard left and right ultrafunctors  as instances of lax/colax algebra morphisms.\end{abstract}
\setcounter{secnumdepth}{2}
\setcounter{tocdepth}{2}
\tableofcontents
\newpage

\section{Introduction}

The ultraproduct construction plays a fundamental role in model theory. Makkai introduced the concept of ultracategories and ultrafunctors \cite{makkai1987stone,makkai88strong}, and used it to show the strong conceptual completeness theorem, which states that for any small pretopos $\ca{T}$, there is an equivalence of categories between $\ca{T}$, and the category of ultrafunctors from $\mathrm{Mod}(\ca{T})$ to  $\cat{Set}$. Here, ultrafunctors are functors that preserve the ultraproduct up to coherence. Early work exploring these notions further includes Marmolejo's 1995 PhD thesis \cite{marmolejo}, which gave a $2$-monadic treatment of Makkai's ultracategories, and Zawadowski \cite{ZawadowskiPhD}, although this early work is not directly related to the results of the present paper. Lurie \cite{lurie2018ultracategories} built on Makkai's work, introducing a much simplified notion of ultracategories, and various notions of functors between these, namely ultrafunctors, right ultrafunctors, and most importantly left ultrafunctors, and gave a different proof of Makkai's conceptual completeness (although Lurie and Makkai's definitions of ultracategories and ultrafunctors differ, they agree in the statement of Makkai's conceptual completeness). Towards this, he showed another important result, which generalises Makkai's conceptual completeness, and Deligne's completeness (the statement that any coherent topos has enough points), and states that for any coherent topos there is an equivalence between $\ca{T}$ and left ultrafunctors from $\mathrm{Mod}(\ca{T})$ to $\cat{Set}$. Left ultrafunctors are a weakening of the definition of ultrafunctors, and unlike ultrafunctors which seem to impose a very restrictive condition on functors (preservation of ultraproduct), are much more abundant in mathematical structures (for example, the $\mathbf{GNS}$ construction in operator algebras is a left ultrafunctor \cite[8.4]{Bundles}). Rewriting ultracategories in a monadic definition was an  early question. In this  direction we mention Clementino and Tholen \cite{CLEMENTINO200313}, who introduce a notion of ultracategories as an instance of a generalised multicategory, in particular a $(\beta \dash \cat{Set})$-multicategory where $\beta$ is the ultrafilter monad, however their notion looked more like generalised ultracategories (which we will talk about in the next paragraph) than ordinary ultracategories. In a blog comment \cite{ShulmanComment}, Shulman postulated that ultracategories are normal colax algebras for a certain pseudo-monad on the category of categories. We prove a corrected version of this conjecture, where we had to make some modifications to the definition suggested by Shulman, to obtain the intended result. We give the definition of the pseudo-monad $T$ in section \ref{pseudo}, and we think of this monad as sending a category to the category of ``formal'' ultraproducts on this category, while the algebra structure in a category recovers the actual ultraproducts (think of any free algebraic structure monad for an oversimplified comparison). In general our approach is also similar to the unrelativised pseudo-monad of the relative ultrafilter monad \cite{tarantino2025ultracategorieskanextensionsrelative}, this pseudo-monad resembles our $T$ except that morphisms are not quotiented by an equivalence relation, the  colax algebras for this pseudo-monad are what the authors call ``weak ultracategories''. Rosolini and Garner independently investigated the same monad of this paper where it was called the ultracompletion pseudo-monad \cite{ultracompletion}, with the direction of their research going into understanding Makkai's ultracategories, and drawing comparison to Lurie's. Indeed, $T$ restricted to the one-object category $\mathbf{1}$ recovers exactly Blass's classical category of ultrafilters  \cite{BlassPhD}. Other work on categorification of ultraproducts include Garners construction of a process $\ca{UF}$, which he used among other things to give a categorical description of ultraproducts of sets \cite{GARNER2020102831}, and also to give a categorical description of the ultraproduct of ultrafilters (called in his case the directed sum of ultrafilters), we will later see that the construction $\ca{UF}$ can be used to understand, for a category $A$ the category $TA$, but without any functoriality or monad structure. Moreover, we should mention that a categorical description of ultraproducts, in this case for all models of coherent logic, was given by Di Liberti \cite{di2022geometry}.

One motivating reason behind introducing this monad will not be clear until our further work which introduces generalised ultracategory in the sense of the author's own work \cite{hamad2025generalisedultracategoriesconceptualcompleteness}, also called virtual ultracategories by Saadia \cite{saadia2025extendingconceptualcompletenessvirtual} and ultraconvergence spaces by Tarantino, Marquès and van Gool \cite{VANGOOL2026108269}, these are relational counterparts to ultracategories in which we don't have an ultraproduct but we have the data of formal ``representables'' at ultraproducts i.e we have for each object $A$ and family of objects $(B_i)_{i \in I}$ a generalised $\la{Hom}$-set, which we denote by $\la{Hom}(A,\int_IB_id\omega)$. It turns out that every pseudo-algebra for the pseudo monad $T$ (which we will now call a strong ultracategory) is in fact a generalised ultracategory, so generalised ultracategories are a generalisation of $T$-pseudo-algebras (strong ultracategories). In Saadia's formulation, the definition is simplified by having a Hom functor $X \times T X \to \cat{Set}$, where $T$ is the present monad restricted to the discrete category on $X$. In a recent paper, Aristote and Tarantino \cite{aristote2026profunctorialalgebras} characterised generalised ultracategories/virtual ultracategories/ ultraconvergence spaces as normalised algebras for an extension of the ultracompletion pseudo-monad $T$ (in their cases called $\beta$) to the category of categories with profunctors and natural transformations $\cat{Prof}$, and strong ultracategories (i.e $T$ pseudo-algebras) are the representables among those.

\subsection{Outline of methodology}
The main result that we want to show is the following:
\begin{theorem*}

Let $\mathrm{Ult}_{L}$ (resp. $\mathrm{Ult}_{R}$, resp. $\mathrm{Ult}$)  denote the bicategory of ultracategories with left ultrafunctors (resp. right ultrafunctors, resp. ultrafunctors) as $1$-cells, and natural transformations (in the appropriate structure-respecting sense)  as $2$-cells, then 
\begin{enumerate}
    \item There is an equivalence of bicategories between $\mathrm{Ult}_{L}$ (resp. $\mathrm{Ult}_{R}$, resp. $\mathrm{Ult}$), and the bicategory of normal colax algebras for the pseudo-monad $T$ as we defined it in \ref{pseudo}, with lax  morphisms as $1$-morphisms (resp. colax  morphisms, resp. pseudo-morphisms) and with the respective $2$-morphisms of those.

    \item This equivalence is realised concretely by a natural bijection between ultrastructures on a category $A$ and  normal colax $T$-algebra structures on $A$, and given two ultracategories $A$ and $B$ by a bijection between  left ultrafunctors, (resp. right ultrafunctors) (resp. ultrafunctors) and left/lax morphisms (resp. right/colax morphisms, resp. pseudo-morphisms) of normal colax algebras,  moreover given two such functors, there is a bijection between natural transformations of left ultrafunctors (resp. right ultrafunctors) (resp. ultrafunctors) between these functors, and the corresponding $2 \dash$ morphisms between the $T$-algebra morphism structures on these functors.
\end{enumerate}

\end{theorem*}

The first few sections give some necessary definitions: in section \ref{Ultra}, we give a comprehensive definition of the $2$-categories of ultracategories in question: ultracategories with left (right) ultrafunctors and natural transformations of left (right) ultrafunctors. Section \ref{colax} provides a definition of various $2$-categories of colax algebras over a $2$-monad.

After that, we give the necessary definitions needed to construct the equivalence. Section \ref{pseudo} gives the definition of the pseudo-monad over the category of locally small categories. This (pseudo)-monad $T$ is defined by sending a category $A$ to the category of tuples of elements of $A$ indexed by a set $I$ together with an ultrafilter $\omega$ on $I$. This closely mirrors the definition given by \cite{ShulmanComment}. For morphisms, one had to reverse the direction of morphisms given by Shulman, together with taking germs of families of maps with respect to the equivalence relation induced by the ultrafilter. This pseudo-monad $T$ was also introduced by Rosolini and Garner \cite{ultracompletion}, where it's called the ultracompletion pseudo-monad, and the authors are interested in establishing a connection between this monad and Makkai's ultracategories \cite{makkai88strong,makkai1987stone}.

In section \ref{First direction}, we construct a functorial assignment sending every ultrastructure on a category $A$ to a $T$-normal colax algebra. An important aspect of this construction is the colax associator, which turns out to be a particular case of the categorical Fubini transform (a part of the ultrastructure), corresponding to the case of a family of ultrafilters $(\gamma_s)$ on sets $(T_s)$ and an ultrafilter $\omega$ on the set $S$, composed with a natural isomorphism; more precisely, the following composition:

\[
\begin{tikzcd}
{\int_{\coprod_{s \in S}T_s}M_{(s,t)}d \int_S \iota_s \gamma_s d \omega } \arrow[rr, "{\Delta_{\omega, \iota \gamma_{\bullet}}}"] &  & {\int_S\int_{\coprod_{s \in S}T_s} M_{(s,t)} d \iota_s \gamma_s d\omega} \arrow[rr, "{\int_S\Delta_{\gamma_s,\iota_s}} d\omega"] &  & \int_S\int_{T_s} M_{(s,t)} d \gamma_s d\omega
\end{tikzcd}\]

One can immediately see two key difficulties: namely one needs to show that this structure is a $T$-colax algebra (we noticed that the colax algebra axioms are ``similar'' to ultracategory axioms, but with the Fubini transform replaced by the colax associator), and we also need to be able to recover the original ultrastructure from the colax associator. 

One important thing we can deduce from this construction is the importance of the colax algebra associator, and a particular family of Fubini transforms involved in its construction: Using the equivalence we show in this paper, pseudo-algebras for this (pseudo)-monad correspond to ultracategories for which these particular families of Fubini transforms are invertible. In section \ref{algebras}, we highlight the fact that all known ultracategories are in fact pseudo-algebras for this pseudo-monad and not just colax algebras (which means that the associator is in fact invertible). Let us explain further, suppose that we have a set $S$, and ultrafilter $\omega$ on $S$, a family of sets $(T_s)_{s \in S}$, together with a family of ultrafilters on each $T_s$, $(\gamma_s)_{s \in S}$; let $\iota_s$ denote the inclusion map of $T_s$ in $\coprod_{s \in S} T_s$, and let $\iota_s \gamma_s$ denote the pushforward of the ultrafilter $\gamma_s$ by the map $\iota_s$, then we can define the categorical Fubini transform \ref{Ultra} $$ \Delta_{\omega,\iota\gamma_\bullet}:\int_{\coprod_{s \in S}T_s}M_{(s,t)}d \int_S \iota_s \gamma_s d \omega  \longrightarrow \int_S\int_{\coprod_{s \in S}T_s} M_{(s,t)} d \iota_s \gamma_s d\omega$$

We see that this map appears in the definition of the colax algebra associator of an ultracategory, so one key step in showing the equivalence, would be to relate a general Fubini transform into one that is of the form above, moreover it's enough that such particular family of Fubini transforms is invertible for the associator to be invertible, so we may define a strong ultracategory ($T$-pseudo algebra) as an ultracategory for which the particular family of Fubini transforms above is invertible.

In section \ref{Second direction}, we construct what will be the inverse of the process in \ref{First direction}. In section \ref{Showing bijection}, we show that these two processes are really inverses of each other, establishing a bijection between these structures on a given category, and then showing that this bijection is functorial. Section \ref{natural transformations} is dedicated to extending this bijection to a full equivalence of $2$-categories for each type of morphism between ultracategories, like left and right ultrafunctors, where these correspond to certain weak morphisms between colax algebras.

\section{Ultracategories} \label{Ultra}

\subsection{Definition of an ultracategory}
Following \cite{lurie2018ultracategories}, we define ultracategories as follows:
\begin{definition}
    
An ultrastructure on a category $A$ consists of the following data:
\begin{enumerate}
    \item For every set $X$, and every ultrafilter $\omega$ on this set, a functor from $A^X$ to $A$ termed the ultraproduct functor, which we denote $\int_X \bullet \, d \omega$.
    \item For every set $X$, every family of ultrafilters on $X$, $(\nu_s)_{s \in S}$, every ultrafilter $\omega$ on $S$, and every family of objects $(M_x)_{x \in X}$, a morphism called the categorical Fubini transform $\Delta_{\omega,\nu_{\bullet}}: \int_X M_x d(\int_S \nu_s d \omega) \rightarrow \int_S (\int_X M_x d \nu_s) d \omega$, which is natural in the family $(M_x)_{x \in X}$.
    \item For every principal ultrafilter $\delta_{x_0}$ on a set $X$, a natural isomorphism $\epsilon_{X,x_0}$ from $\int_X M_x d \delta_{x_0}$ to $M_{x_0}$.
\end{enumerate}
This data is required to satisfy the following axioms: 

\begin{enumerate}[label=\Alph*]
\item Suppose we have a family of ultrafilters $(\nu_s)_{s \in S}$, then the map \newline $\Delta_{\delta_{s_0},\nu_{\bullet}}: \, \int_X M_x d \int_S \nu _{s} d \delta_{s_{0}} \rightarrow \int_S \int_X M_x d \nu_s d \delta_{s_{0}}$ is the inverse of the map $\epsilon_{S,s_0}$ from  $\int_S \int_X M_x d \nu_s d \delta_{s_{0}}$ to $\int_X M_x d \nu_{s_0}$.

\item Suppose we have a monomorphism of sets $f: \, X \rightarrow Y $, and an ultrafilter $\omega$ on $X$, then the categorical Fubini transform from  $\int_Y M_y df \omega=\int_Y M_y d \int_X \delta_{f(x)} d \omega$ to $\int_X \int_Y M_y d \delta_{f(x)} d \omega$ is an isomorphism. 

Here $f\omega$ is the pushforward of $\omega$ by the map $f$ defined by $B \in f\omega \iff f^{-1}(B) \in \omega$.
\item Suppose we have a set $R$ and an ultrafilter $\lambda$ on $R$, and suppose we have $(\omega_r)_{r \in R}$ a family of ultrafilters on a set $S$, and $(\nu_s)_{s \in S}$ is a family of ultrafilters on some set $T$, then if we define 

$\rho=\int_R (\int_S \nu_s d\omega_r)d\lambda =\int_S \nu_s d (\int_R \omega_r d \lambda)$, here $\int_S \nu_s d \omega$ is defined by $$B \in \int_S \nu_s d \omega \iff \{s \in S: \ B \in \nu_s \} \in \omega$$ then the following diagram commutes:
\[
\begin{tikzcd}
\int_{T} M_t d \rho \arrow[rrr, "{\Delta_{\lambda,\int_S \nu_{s}d\omega_{\bullet}}}"] \arrow[dd, "{\Delta_{\int_R \omega_r d \lambda,\nu_{\bullet}}}"] &  &  & \int_{R} (\int_{T} M_t d \int_S \nu_s d \omega_r) d \lambda \arrow[dd, "{\int_R \Delta_{\omega_r,\nu_\bullet}d \lambda}"] \\
& & & \\
\int_S \int_T M_t d \nu_s d \int_R \omega_r d \lambda \arrow[rrr, "{\Delta_{\lambda,\omega_{\bullet}}}"]                                               &  &  & \int_{R} (\int_S (\int_T M_t d \nu_s)d \omega_r) d \lambda                                                             
\end{tikzcd}
\]

\end{enumerate}
\end{definition}

Now we define an ultracategory to be a category with an ultrastructure.

Given an ultracategory, there is an important class of maps that is not part of the ultracategory definition but still plays an important role, these are called the ultraproduct diagonal maps $\Delta_{\omega,f}$, Suppose that we have an ultrafilter $\omega$ on $J$ and a map of sets $f$ from  $J$ to $I$, then we define $\Delta_{\omega,f}$ from $\int_I M_i df\omega$ to $\int_J M_{f(j)} d \omega$ by taking the following composition:

\[\begin{tikzcd}
	{\int_I M_i d f\omega = \int_I M_i d \int_J \delta_{f(j)} d\omega} &&&& {\int_J\int_IM_i d \delta_{f(j)}d\omega} &&&& {\int_J M_{f(j)}d\omega}
	\arrow["{\Delta_{\omega,\delta_{f \bullet}}}", from=1-1, to=1-5]
	\arrow["{\int_J \epsilon_{I, \delta_{f(j)}} d\omega}", from=1-5, to=1-9]
\end{tikzcd}\]

\subsection{Left and right ultrafunctors}
Let $\mathsf{M}$ and $\mathsf{N}$ be two ultracategories. A left ultrafunctor from $\mathsf{M}$ to $\mathsf{N}$ is a functor equipped with a left ultrastructure.

\begin{definition}
   A left ultrastructure on a functor $F$ consists of the following: for every ultrafilter $\omega$ on a set $X$, and for every family of objects $(M_x)_{x \in X}$ in $\mathsf{M}$, there is a  morphism $F(\int_X M_x d\omega) \to \int_X F(M_x)d \omega$, which we collectively denote by $\sigma_\omega$, although it depends on the family $(M_x)_{x \in X}$, as well as on $\omega$, since the family will always be clear from context. 

    These morphisms satisfy the following axioms:
    \begin{enumerate}
        \item  For every family of morphisms $(\psi_x)_{x \in X}$ from $M_x$ to $N_x$ in $\mathsf{M}$, the following diagram commutes:
        \[
            \begin{tikzcd}
                F(\int_X M_x d \omega) \arrow[rr, "\sigma_{\omega}"] \arrow[dd, "F(\int_X \psi_x  d \omega)"] &  & \int_XF(M_x)d \omega \arrow[dd, "\int_XF(\psi_x) d \omega"] \\
                & & \\
                F(\int_X N_x d \omega) \arrow[rr, "\sigma_{\omega}"] &  & \int_XF(N_x)d\omega
            \end{tikzcd}
        \]

        \item For every set $X$ and every $x_0 \in X$, the following diagram commutes:
        
\[\begin{tikzcd}[ampersand replacement=\&]
	{F(\int_X M_x d \delta_{x_{0}})} \&\& {\int_X F(M_x) d \delta_{{x_0}}} \\
	\& {F(M_{x_0})}
	\arrow["{\sigma_{\delta_{x_{0}}}}", from=1-1, to=1-3]
	\arrow["{{F(\epsilon_{X,x_{0}})}}"'{pos=0.3}, from=1-1, to=2-2]
	\arrow["{{\epsilon_{X,x_{0}}}}", from=1-3, to=2-2]
\end{tikzcd}
\]

        \item Given sets $S$ and $T$, an ultrafilter $\omega$ on $S$,  and a family of ultrafilters $(\nu_s)_{s \in S}$ on $T$, the following diagram commutes:
        \[
            \begin{tikzcd}
                F(\int_T M_t d (\int_S \nu_s d \omega)) \arrow[rrrr, "\sigma_{\int_S\nu_s d \omega}"] \arrow[dd, "{F(\Delta_{\omega,\nu_{\bullet}})}"] & & & & \int_T F(M_t) d \int_S \nu_s d\omega \arrow[dd, "{\Delta_{\omega,\nu_{\bullet}}}"] \\
                & & & & \\
                F(\int_S (\int_T M_t d\nu_s) d\omega) \arrow[rr, "\sigma_{\omega}"] & & \int_SF(\int_TM_td\nu_s)d \omega \arrow[rr, "\int_S \sigma_{\nu_s} d \omega"] & & \int_S \int_T F(M_t) d \nu_s d \omega
            \end{tikzcd}
        \]

    \end{enumerate}
\end{definition}

\begin{note*}
    \normalfont
    The dual notion is a right ultrafunctor, in which the comparison maps go in the other direction; we omit writing the axioms which can be found in \cite[8.1]{lurie2018ultracategories}.
\end{note*}

\begin{note*}
    An ultrafunctor is a left (or right) ultrafunctor for which the morphisms $\sigma_{\omega}$ are invertible.
\end{note*}

\subsection{Natural transformations of ultrafunctors}
Suppose that $\mathsf{M}$ and $\mathsf{N}$ are two ultracategories, and let $F$, $G$ be left ultrafunctors between them. A natural transformation of left ultrafunctors from $F$ to $G$ is a natural transformation $\phi$ satisfying the additional condition: 

For every family $(M_i)$ of objects in $\mathsf{M}$ and for every ultrafilter $\omega$ on $I$, the following diagram commutes (here $\sigma$ and $\sigma^{'}$ denote the left ultrastructure maps of $F$ and $G$ respectively):
\[
\begin{tikzcd}[ampersand replacement=\&]
	{F(\int_I M_i d\omega)} \&\&\& {\int_I F(M_i)d\omega} \\
	\\
	{G(\int_I M_i d\omega)} \&\&\& {\int_I G(M_i)d\omega}
	\arrow["{{\sigma_{\omega}}}"', from=1-1, to=1-4]
	\arrow["{{\phi_{\int_I M_i d\omega}}}", from=1-1, to=3-1]
	\arrow["{{\int_I \phi_{M_i}d\omega}}", from=1-4, to=3-4]
	\arrow["{{\sigma^{'}_{\omega}}}", from=3-1, to=3-4]
\end{tikzcd}
\]
A natural transformation of right ultrafunctors is defined similarly.

 \section{Colax algebra over a monad} \label{colax} 

 We are going to give the definition of the category of algebras for strict $2$-monads, since we can reduce our study to those. It is still possible to work with pseudo-monads and study algebras of those, for a reference on these see \cite{lack2000coherent} (notice that the author only works with pseudo-algebras \cite[4.2]{lack2000coherent}, but it's not hard to weaken these to obtain various categories of lax/colax algebras).
 
 Following \cite[section 3]{penon2009representable} and \cite[1.1]{CTGDC}, we define a colax algebra over a monad as follows (the paper \cite{penon2009representable} calls such  algebras algèbres laxes (lax algebras), while they are called pseudo-algèbres à droite (right pseudo-algebras) in \cite[1.1]{CTGDC}, we will call them colax algebras): Suppose that we have a $2$-monad $(M,\eta,\mu)$ over some $2$-category $\ca{C}$, then a colax algebra is an object $C \in \ca{C}$ equipped with:

 \begin{itemize}

\item  $1$-morphism $m$:  $MC \xrightarrow{} C$, which is called the \emph{algebra action}.
\item  $2$-morphism $i$: 
\[\begin{tikzcd}[ampersand replacement=\&,cramped]
	C \&\& MC \\
	\\
	\&\& C
	\arrow["{\eta_C}"', from=1-1, to=1-3]
	\arrow[""{name=0, anchor=center, inner sep=0}, "{\la{Id}}"', from=1-1, to=3-3]
	\arrow["m"', from=1-3, to=3-3]
	\arrow["i"', shorten >=5pt, Rightarrow, from=1-3, to=0]
\end{tikzcd}\] which we call the \emph{unitor}.
\item   $2$-morphism  $a$: 
\[\begin{tikzcd}[ampersand replacement=\&]
	{M^2C} \&\&\& MC \\
	\\
	MC \&\&\& C
	\arrow["Mm", from=1-1, to=1-4]
	\arrow["{\mu_{C}}", from=1-1, to=3-1]
	\arrow["m"', from=1-4, to=3-4]
	\arrow["a", shorten >=16pt, Rightarrow, from=3-1, to=1-4]
	\arrow["m", from=3-1, to=3-4]
\end{tikzcd}\]
which we call the \emph{associator}.
 \end{itemize}
All are required to satisfy the following:

\adjustbox{max width=\textwidth}{

\begin{tikzcd}[ampersand replacement=\&]
	MC \&\& {M^2C} \&\& MC \&\&\& {=} \& MC \&\& {M^2C} \&\& MC \&\& C \& {=} \& MC \&\&\& C \\
	\&\& C \&\& MC \&\& C \&\&\&\&\&\& MC
	\arrow["{\eta_{MC}}", from=1-1, to=1-3]
	\arrow["m"', from=1-1, to=2-3]
	\arrow["{\mu_C}", from=1-3, to=1-5]
	\arrow["Mm"{description}, from=1-3, to=2-5]
	\arrow["a"', Rightarrow, nfold, from=1-5, to=2-5]
	\arrow["m", from=1-5, to=2-7]
	\arrow["{M\eta_C}", from=1-9, to=1-11]
	\arrow[""{name=0, anchor=center, inner sep=0}, "{\la{Id}}"', curve={height=18pt}, from=1-9, to=2-13]
	\arrow["{\mu_C}", from=1-11, to=1-13]
	\arrow["Mm"', from=1-11, to=2-13]
	\arrow["m", from=1-13, to=1-15]
	\arrow["a"', Rightarrow, nfold, from=1-13, to=2-13]
	\arrow["m", from=1-17, to=1-20]
	\arrow["{\eta_C}"', from=2-3, to=2-5]
	\arrow[""{name=1, anchor=center, inner sep=0}, "{\la{Id}}"', curve={height=30pt}, from=2-3, to=2-7]
	\arrow["m", from=2-5, to=2-7]
	\arrow["m"', from=2-13, to=1-15]
	\arrow["Mi"', between={0}{0.85}, Rightarrow, nfold, from=1-11, to=0]
	\arrow["i"{description}, between={0}{0.8}, Rightarrow, nfold, from=2-5, to=1]
\end{tikzcd}
}
and

\adjustbox{max width =\textwidth}{
\begin{tikzcd}
	{M^3C} && {M^2C} && MC && C & {=} & {M^3C} && {M^2C} && MC && C \\
	&& {M^2C} && MC && {} &&&& {M^2C} && MC \\
	&&&&&&&&&&&& MC
	\arrow["{M\mu_C}", from=1-1, to=1-3]
	\arrow["{\mu_C}", from=1-3, to=1-5]
	\arrow["m", from=1-5, to=1-7]
	\arrow["{M m}", from=1-3, to=2-5]
	\arrow["m", from=2-5, to=1-7]
	\arrow["{M^2m}"', from=1-1, to=2-3]
	\arrow["Mm"', from=2-3, to=2-5]
	\arrow["Ma"', Rightarrow, from=1-3, to=2-3]
	\arrow["a"', Rightarrow, from=1-5, to=2-5]
	\arrow["{\mu_{MC}}", from=1-9, to=1-11]
	\arrow["{\mu_C}", from=1-11, to=1-13]
	\arrow["m", from=1-13, to=1-15]
	\arrow["Mm"{description}, from=1-11, to=2-13]
	\arrow["m"{description}, from=2-13, to=1-15]
	\arrow["a", Rightarrow, from=1-13, to=2-13]
	\arrow["{M^2m}", from=1-9, to=2-11]
	\arrow["{\mu_C}", from=2-11, to=2-13]
	\arrow["{M m}"{description}, from=2-11, to=3-13]
	\arrow["m"{description}, from=3-13, to=1-15]
	\arrow["a", Rightarrow, from=2-13, to=3-13]
\end{tikzcd}
}

A colax algebra over a monad is called \emph{normal} if the unitor $i$ is an isomorphism.
\subsection{The 2-category of colax algebras over a monad}

Suppose that we are given two colax algebras $(A,m_1,a_1,i_1)$ and $(B,m_2,a_2,i_2)$ over a monad $M$ on a $2$-category $\ca{C}$, we define a lax (or left)  morphism of colax algebras between $A$ and $B$  to be a morphism in $\ca{C}$ between $A$ and $B$ together with a 2-morphism $\phi$: $f \circ m_1 \rightarrow m_2 \circ Mf$ satisfying the following :

\adjustbox{max width = \textwidth}{
\begin{tikzcd}
	A && MA && A &&& {=} & A &&& MA &&& A &&& B \\
	&& B && MB && B
	\arrow["{{\eta_A}}", from=1-1, to=1-3]
	\arrow["f"', from=1-1, to=2-3]
	\arrow["{m_1}", from=1-3, to=1-5]
	\arrow["Mf"', from=1-3, to=2-5]
	\arrow["\phi"{description}, Rightarrow, from=1-5, to=2-5]
	\arrow["f", from=1-5, to=2-7]
	\arrow["{{\eta_A}}", from=1-9, to=1-12]
	\arrow[""{name=0, anchor=center, inner sep=0}, curve={height=30pt}, from=1-9, to=1-15]
	\arrow["{m_1}", from=1-12, to=1-15]
	\arrow["f"', from=1-15, to=1-18]
	\arrow["{{\eta_B}}"', from=2-3, to=2-5]
	\arrow[""{name=1, anchor=center, inner sep=0}, curve={height=30pt}, from=2-3, to=2-7]
	\arrow["{m_2}"', from=2-5, to=2-7]
	\arrow["{i_1}"', between={0}{0.8}, Rightarrow, from=1-12, to=0]
	\arrow["{i_2}"{description}, between={0}{0.8}, Rightarrow, from=2-5, to=1]
\end{tikzcd}
}

\adjustbox{max width =\textwidth}{

\begin{tikzcd}
	{M^2A} && MA && A && B && {M^2A} &&& MA && A && B \\
	& {M^2 B} && MB &&&& {=} &&&& MA && MB \\
	&&& MB &&&&&&& {M^2B}
	\arrow["{{\mu_A}}", from=1-1, to=1-3]
	\arrow["{{M^2f}}"', from=1-1, to=2-2]
	\arrow[""{name=0, anchor=center, inner sep=0}, "{m_1}", from=1-3, to=1-5]
	\arrow["Mf"{description}, from=1-3, to=2-4]
	\arrow["f", from=1-5, to=1-7]
	\arrow["{{\mu_{A}}}", from=1-9, to=1-12]
	\arrow["{Mm_1}"{description}, from=1-9, to=2-12]
	\arrow["{{M^2f}}"{description}, from=1-9, to=3-11]
	\arrow["{m_1}", from=1-12, to=1-14]
	\arrow["a", Rightarrow, from=1-12, to=2-12]
	\arrow["f", from=1-14, to=1-16]
	\arrow["\phi", Rightarrow, from=1-14, to=2-14]
	\arrow["{{\mu_B}}"', from=2-2, to=2-4]
	\arrow["{{M m_2}}"', from=2-2, to=3-4]
	\arrow["{m_2}"{description, pos=0.4}, shift left, from=2-4, to=1-7]
	\arrow["{a_2}", Rightarrow, from=2-4, to=3-4]
	\arrow["{m_1}"{description}, from=2-12, to=1-14]
	\arrow["Mf"{description}, from=2-12, to=2-14]
	\arrow["{M\phi}"{description}, Rightarrow, from=2-12, to=3-11]
	\arrow["{m_2}"', from=2-14, to=1-16]
	\arrow["{ m_2}"', from=3-4, to=1-7]
	\arrow["{M m_2}"{description}, from=3-11, to=2-14]
	\arrow["\phi", between={0.2}{1}, Rightarrow, from=0, to=2-4]
\end{tikzcd}
}

\begin{note*} \normalfont In the paper of Burroni \cite[1.2]{CTGDC}, natural transformations $\phi$ going the other direction are called morphismes à gauche (left morphisms), and the ones we defined here are called morphismes à droite (right morphisms), but we decided to call them left morphisms, since they are going to correspond to left ultrafunctors between ultracategories in the special case we are treating, we also call them lax morphisms following \cite[page 3]{BLACKWELL19891}. The adjective lax has nothing to do with the fact that we have colax algebras over a monad.
\end{note*}

Alternatively, one can define a colax (right)  morphism of colax algebras (or left morphisms of colax algebras if one follows the terminology of \cite{CTGDC}) as a morphism $f$ in $\ca{C}$ between $A$ and $B$, together with a $2$-morphism $\phi:m_2 \circ Mf \rightarrow f \circ m_1 $, satisfying the following:

\adjustbox{max width = \textwidth}{

\begin{tikzcd}
	A && MA &&& MB &&& {=} & A &&& MA &&& MB && B \\
	&& B &&&&& B &&&&&& A
	\arrow["{{\eta_A}}", from=1-1, to=1-3]
	\arrow["f"', from=1-1, to=2-3]
	\arrow["Mf", from=1-3, to=1-6]
	\arrow["{m_2}", from=1-6, to=2-8]
	\arrow["{{\eta_A}}", from=1-10, to=1-13]
	\arrow[""{name=0, anchor=center, inner sep=0}, "{{\la{Id}}}"{description}, from=1-10, to=2-14]
	\arrow["Mf", from=1-13, to=1-16]
	\arrow["{m_1}", from=1-13, to=2-14]
	\arrow["{{m_2}}", from=1-16, to=1-18]
	\arrow["\phi"{description}, Rightarrow, from=1-16, to=2-14]
	\arrow["{{\eta_B}}", from=2-3, to=1-6]
	\arrow[""{name=1, anchor=center, inner sep=0}, "{{\la{Id}}}"', from=2-3, to=2-8]
	\arrow["f"{description}, from=2-14, to=1-18]
	\arrow["{{i_2}}", between={0}{0.8}, Rightarrow, from=1-6, to=1]
	\arrow["{i_1}"{description}, between={0}{0.8}, Rightarrow, from=1-13, to=0]
\end{tikzcd}

}

\adjustbox{max width = \textwidth}{
\begin{tikzcd}
	&& MA \\
	{M^2A} && MA && A && B && {M^2A} &&&&& A && B \\
	&&&&&&& {=} &&&& MA && MB \\
	{M^2 B} &&& MB &&&&&&& {M^2B} &&& MB
	\arrow["{m_1}"{description}, from=1-3, to=2-5]
	\arrow["{Mm_1}"{description}, from=2-1, to=1-3]
	\arrow["{{\mu_{A}}}", from=2-1, to=2-3]
	\arrow["{{M^2f}}"', from=2-1, to=4-1]
	\arrow["a"', Rightarrow, from=2-3, to=1-3]
	\arrow[""{name=0, anchor=center, inner sep=0}, "{m_1}", from=2-3, to=2-5]
	\arrow["Mf"{description}, from=2-3, to=4-4]
	\arrow["f", from=2-5, to=2-7]
	\arrow["{Mm_1}"{description}, from=2-9, to=3-12]
	\arrow["{M^2f}"{description}, from=2-9, to=4-11]
	\arrow["f", from=2-14, to=2-16]
	\arrow["{m_1}"{description}, from=3-12, to=2-14]
	\arrow["Mf"{description}, from=3-12, to=3-14]
	\arrow["\phi"', Rightarrow, from=3-14, to=2-14]
	\arrow["{m_2}"{description}, from=3-14, to=2-16]
	\arrow["{{\mu_B}}"{description}, from=4-1, to=4-4]
	\arrow["{m_2}"'{pos=0.4}, shift left, from=4-4, to=2-7]
	\arrow["{M\phi}", shift right, Rightarrow, from=4-11, to=3-12]
	\arrow["{M m_2}"{description}, from=4-11, to=3-14]
	\arrow["{{\mu_B}}"', from=4-11, to=4-14]
	\arrow["{m_2}"', from=4-14, to=2-16]
	\arrow["{Ma_2}"', Rightarrow, from=4-14, to=3-14]
	\arrow["\phi", Rightarrow, from=4-4, to=0]
\end{tikzcd}
}
\begin{note*} \normalfont
    In our case, we defined lax (left) and colax (right) morphisms between colax algebras over a strict $2$-monad, similarly one can define these concepts for lax algebras over a pseudo-monad.
\end{note*}

\noindent Of course, if we ask the $2$-morphism $\phi$ in either of the two definitions above to be invertible, we get what we are going to call a pseudo-morphism of colax algebras.

Now we define the $2$-morphisms between colax morphisms, suppose that we have $f$, $g$ two left morphisms of colax algebras, then we define a $2$-morphism in that $2$-category to be a $2$-morphism $\psi$ in the category $\ca{C}$ between $f$ and $g$ (to be more precise, between the underlying $1$-morphisms of $f$ and $g$ in the category $\ca{C}$, which we also denote by $f$ and $g$) satisfying the following coherence condition which can be expressed by the commutativity of the following diagram:

\[\begin{tikzcd}
	{m_2 \circ Mf} &&& {m_2 \circ Mg} \\
	\\
	{f \circ m_1} &&& {g \circ m_1}
	\arrow["{{m_2 \circ M\psi}}"', from=1-1, to=1-4]
	\arrow["\phi"', from=3-1, to=1-1]
	\arrow["{{\psi \circ m_1}}", from=3-1, to=3-4]
	\arrow["{{\phi^{'}}}", from=3-4, to=1-4]
\end{tikzcd}\]

\noindent Alternatively, if we have two right morphisms of colax algebras, we ask the commutativity of the following diagram :

\[\begin{tikzcd}
	{m_2 \circ Mf} &&& {m_2 \circ Mg} \\
	\\
	{f \circ m_1} &&& {g \circ m_1}
	\arrow["{{m_2 \circ M\psi}}"', from=1-1, to=1-4]
	\arrow["\phi", from=1-1, to=3-1]
	\arrow["{{\phi^{'}}}"', from=1-4, to=3-4]
	\arrow["{{\psi \circ m_1}}", from=3-1, to=3-4]
\end{tikzcd}\]

\subsection{The Monad Definition} \label{pseudo}
We define a pseudo-monad $T$ on the $2$-category of locally small categories $\cat{CAT}$ as follows: Suppose $U$ is a locally small category. The category $TU$ has for objects $(X,\omega,(M_x)_{x \in X})$ where $X$ is a set, $(M_x)$ is a family of objects of $U$ indexed by $X$, and $\omega$ is an ultrafilter on $X$.

A morphism in the category $TU$ has the following data: Suppose that $(X,\omega,(M_x)_{x \in X})$ and $(Y,\nu,(N_y)_{y \in Y})$ are objects, then a morphism is a pair $(f,(\sigma_{y}))$ which denotes the equivalence class of set functions $f$ from $ Y^{'} \subseteq Y $ to $X$, that satisfies the following:  $Y^{'} \in \nu$, and if we call $\nu^{'}$ the restriction of $\nu$ to $Y^{'}$, we have $f \nu^{'}=\omega$, subject to the following equivalence relation: $f$ is equivalent to $g$, if $f$ and $g$ agree on some $M\subseteq Y^{'} \bigcap Y^{''}$ such that $M \in \nu$. 

The family $(\sigma_{y})_{y \in Y^{'}}$ (here $Y^{'}$ could be any subset of $Y$ such that $Y^{'} \in \nu$) is a family of morphisms from $M_{f(y)}$ to $N_y$ modulo the equivalence relation defined by identifying two such families  if they agree on a subset $M^{'} \in \nu$. In other words $$\la{Hom}(( X,\omega,(M_x)_{x \in X}),(Y,\nu,(N_y)_{y \in Y})) = \{ (f: Y^{'} \xrightarrow{} X, (\sigma_y: M_{f(y)} \xrightarrow{} N_y)), Y^{'} \in \nu, \ f\nu^{'} = \omega \}/\sim  $$  \noindent where $\sim$ is the equivalence relation $(f,(\sigma_{y})) \sim (f^{'},(\sigma^{'}_{y}))$ if there exists a subset $Y^{''} \in \nu$ such that   $f|_{Y^{''}} = f^{'}|_{Y^{''}}$ and for every $y \in Y^{''}$ $\sigma_{y} =\sigma^{'}_{y}$.

The category $TU$ can be regarded as the category of ``formal ultraproducts'' of elements of $U$, and one can expect the adequate notion of algebra functor to send this formal ultraproduct to the actual one. 

\begin{note*} \normalfont
    Another way to describe the category $TU$ for a category $U$ is through Garner process called $\ca{UF}$ \cite[Definition 21]{GARNER2020102831}, namely given an extensive  category $A$, we may consider the category of ultrafilters on this category called $\ca{UF}(A)$. Now consider the category $Fam(A^{op})$, whose objects are families of objects of $A$ indexed by sets,  and whose morphisms between $(M_i)_{i \in I}$ and $(N_j)_{j \in J}$ are pairs $(f, (\sigma_{i}))$, where $f: I \to J$ and $\sigma_i:  N_{f(i)} \to M_i$ in $A$, this category is extensive \cite{CARBONI1993145}, so we may consider $\ca{UF}(Fam(A^{op}))$, and it's not hard to check that $TA \simeq \ca{UF}(Fam(A^{op}))^{op}$. Unfortunately, we don't see how to use this observation to make $T$ into a monad, and so we didn't make any direct use of it.
\end{note*}

\begin{note*} \normalfont A simpler version of the definition above was first proposed by Shulman \cite{ShulmanComment}. The exact definition above was also given independently, by Rosolini and Garner (unpublished, but already presented see \cite{ultracompletion}) where it's called \emph{the ultracompletion pseudo-monad}, they noted that $TU$ is always an ultracategory, not only in the sense of Lurie as our work may suggest, but also in the sense of Makkai.
\end{note*}

Now suppose that we have a functor $\ca{F}$ from $U$ to $V$ then we get a functor from $TU$ to $TV$ defined by $T(\ca{F})(X,\omega,(M_x)_{x \in X})=(X,\omega,\ca{F}(M_x)) $ and on morphisms by $T(\ca{F})(f,(\sigma_{y}))=(f,(\ca{F}(\sigma_y)))$. Finally, suppose that we have a natural transformation $\alpha $ between two functors $\ca{F}$ and $\ca{G}$ from category $A$ to category $B$:
\[\begin{tikzcd}
	A && B
	\arrow[""{name=0, anchor=center, inner sep=0}, "F", from=1-1, to=1-3]
	\arrow[""{name=1, anchor=center, inner sep=0}, "G"', curve={height=18pt}, from=1-1, to=1-3]
	\arrow["\alpha", shorten <=2pt, shorten >=2pt, Rightarrow, from=0, to=1]
\end{tikzcd}\]
We define the components of a natural transformation between $TF$ and $TG$,  $T\alpha$, as follows :

\[\begin{tikzcd}
	{(X, \omega,\mathcal{F}(M_x)) } \\
	\\
	\\
	{(X, \omega,\mathcal{G}(M_x)) }
	\arrow["{(\la{Id},(\alpha_{M_x}))}", from=1-1, to=4-1]
\end{tikzcd}\]

\begin{note*} \normalfont We leave it to the reader showing that this construction of $T$ is a strict $2$-functor on the strict $2$-category of locally small categories with functors and natural transformations.
\end{note*}

Now we define the monad multiplication $\mu$ as follows: Suppose that we have a category $A$,  an object in $T^2A$ looks like this:
$$(S,\xi,( (T_s)_{s \in S}, (\nu_s)_{s \in S}, (M_{(s,t)})_{(s,t) \in \coprod_{s \in S} T_s}))$$

Here $S$ is a set, $\xi$ an ultrafilter on $S$, $(T_s)_{s \in S}$ is a family of sets indexed by $S$, $(\nu_s)$ is a family of ultrafilters such that $\nu_s$ is an ultrafilter on $T_s$, $(M_{(s,t)})$ is a family of objects of $A$ indexed by $\coprod_{s \in S}T_s$.

For every ultrafilter $\nu_s$ on $T_s$ one can construct an ultrafilter $\iota_s \nu_s$ on $\coprod_{s \in S}T_s $, by pushforward by the map $\iota_s$ of $\nu_s$, and hence we may define the ultrafilter $\int_{S}\iota_s \nu_s d\xi$ on $\coprod_{s \in S}T_s$, which is defined by $A \in \int_{S}\iota_s \nu_s d\xi$ iff $$\{ s \in S | \  A \bigcap T_s \in \nu_s\} \in \xi$$

We send this object by $\mu_A$ to $(\coprod_{s \in S} T_s,\int_S \iota_s \nu_s d \xi,(M_{(s,t)})_{(s,t) \in \coprod_{s \in S}T_s})$, and we define the monad unit $\eta_A$ by sending $X \in A$ to $(*,*,(X))$ where $*$ denotes at the same time the one-object set and the unique ultrafilter on it.

\paragraph{Digression}
The construction above is not a strict $2$-monad on the $2$-category of categories with functors, but rather a pseudo-monad, and hence the definition is incomplete since there are invertible $3$-morphisms that show up as part of the structure when checking monadicity. However, this can be omitted by virtue of the next construction: If we want the construction above to be a strict $2$-monad, we can use the following trick: Instead of defining the object of the $TA$ to be triples $(X,\omega,(M_x)_{x \in X})$ such that $X$ is a set and $\omega$ is an ultrafilter on $X$, we define a strict version of $T$, which we call $T^{'}$:

We take objects of $T^{'} A$ to be  $(X,\omega,(M_x)_{x \in X})$ such that $X$ is an ordinal number, and $\omega$ is an ultrafilter on $X$, the definition of morphisms is the same, but now with this definition, we can define the monad multiplication  by sending:

$$(S,\omega,( (T_s)_{s \in S}, (\nu_s)_{s \in S}, (M_{(s,t)})_{(s,t) \in \coprod_{s \in S} T_s}))$$ to

$$(\Omega,\int_S \iota_s \nu_s d \omega,(M_{(s,t)})_{(s,t) \in \coprod_{s \in S}T_s})$$ where $\Omega$ is the unique ordinal which is order equivalent to $\coprod_{s \in S}T_s$ when $\coprod_{s \in S} T_s$ is equipped with the lexicographic order, and $\int_S \iota_s \nu_s d \omega$ is in reality the image of $\int_S \iota_s \nu_s d \omega$ by the unique order preserving isomorphism, the monad unit definition remains the same (we can choose the one-point set in our previous definition to be the ordinal $\mathbf{1})$. The definition of monad multiplication $\mu$ and unit $\eta$ remain the same, and we are going to be using the same symbol for them.

We can easily see that $T$ and $T^{'}$ are equivalent (and we are going to use an equivalence such that for any category $A$, the equivalence between $T^{'}A$ and $TA$ will send families indexed by ordinals to themselves). This would imply that they have the same category of (colax) algebras, which are our objects of interest. We want to insist on two facts, the first is that $T'$ is an ad-hoc tool intended to simplify our work and obtain a strict $2$-monad, and the second is that, as we are going to see, even if we are given a $T'$-algebra structure on a certain category, we are still interested in $T$-algebra structure as this is the way to obtain ultraproducts for all sets not merely ordinals: 

Since we don't want to give the exact definition of a $T$-colax algebra (since it involves these invertible $3$-morphisms that we don't want to keep track of), we give the following  equivalent definition:

\begin{definition}
    We say that a functor $m: TA \rightarrow A$ has the structure of a $T$-normal colax algebra if there exists a $T^{'}$-normal colax algebra $m^{'}:T^{'}A \rightarrow A$ (which comes with a $2$-morphism $a$ and an invertible $2$-morphism $i$) such that the following diagram commutes :
\[\begin{tikzcd}
	{T^{'}A} && A \\
	TA
	\arrow["{m^{'}}", from=1-1, to=1-3]
	\arrow["\simeq"', from=1-1, to=2-1]
	\arrow["m"', from=2-1, to=1-3]
\end{tikzcd}\]
\end{definition}

Similarly, we can define lax, colax and pseudo-morphisms of normal colax algebras and the corresponding $2$-morphisms (for example a lax morphism between $A$ and $B$ is a functor $f$ together with natural transformation $\phi$ : $f \circ m_1 \rightarrow m_2 \circ Tf$ such that restriction to ordinals would define a lax morphism of colax algebras between $(A,m^{'}_1)$ and $(B,m^{'}_2)$). Although this is not the exact definition of colax algebras for a pseudo-monad, and their lax, colax and pseudo-morphisms, with their corresponding natural transformations, this definition gives an equivalent $2$-category of normal colax algebras.

 \subsection{$T^{'}$ is a monad}
 In what follows, we are going to show that $T^{'}$ is a strict $2$-monad, by checking the commutativity of the following diagrams :

\[\begin{tikzcd}
	{{T^{'}}^{3}} && {{T^{'}}^{2}} && {{T^{'}}} && {{T^{'}}^{2}} \\
	\\
	{{T^{'}}^{2}} && {{T^{'}}} && {{T^{'}}^{2}} && {{T^{'}}}
	\arrow["{{T^{'}\mu}}", from=1-1, to=1-3]
	\arrow["{{\mu T^{'}}}"', from=1-1, to=3-1]
	\arrow["{{\mu}}", from=1-3, to=3-3]
	\arrow["{{T^{'} \eta}}", from=1-5, to=1-7]
	\arrow["{{\eta T^{'}}}"', from=1-5, to=3-5]
	\arrow["{{\la{Id}}}", from=1-5, to=3-7]
	\arrow["{{\mu}}", from=1-7, to=3-7]
	\arrow["{{\mu}}"', from=3-1, to=3-3]
	\arrow["{{\mu}}"', from=3-5, to=3-7]
\end{tikzcd}\]
For an object  $(X,\omega,(Y_x,\lambda_x ,(Z_{(x,y)},\xi_{(x,y)},M_{(x,y,z)}))) \in {T^{'}}^3A$, applying both directions in the first diagram will get us:

\adjustbox{max width = \textwidth}{
\begin{tikzcd}[column sep=20,row sep=40] 
	{(X,\omega,(Y_x,\lambda_x ,(Z_{(x,y)},\xi_{(x,y)},(M_{(x,y,z)}))))} &&& {(\coprod_{x \in X}Y_x,\int_{X}\lambda_{x}d\omega,(Z_{x,y},\xi_{(x,y)},(M_{(x,y,z)})))} \\
	\\
	{(X,\omega,(\coprod_{y\in Y_x}Z_{(x,y)},\int_{Y_x} \iota_{y} \xi_{(x,y)} d \lambda_x,(M_{(x,y,z)})))} &&& {(\coprod_{x \in X}\coprod_{y \in Y_x}Z_{(x,y)},\int_X \iota_x \int_{Y_x}\iota_{y} \xi_{(x,y)} d \lambda_xd\omega,(M_{(x,y,z)})) = (\coprod_{(x,y)\in\coprod_{x \in X}Y_x}Z_{(x,y) }, \int_{\coprod_{x \in X}Y_x}\iota_{(x,y)}\xi_{(x,y)}d\int_{X}\lambda_{x}d\omega,(M_{(x,y,z)}))}
	\arrow[maps to, from=1-1, to=3-1]
	\arrow[maps to, from=3-1, to=3-4]
	\arrow[maps to, from=1-1, to=1-4]
	\arrow[maps to, from=1-4, to=3-4]
\end{tikzcd}

}

Of course, the last equality requires us to take the ordinal which is order isomorphic to  \newline $\coprod_{(x,y)\in\coprod_{x \in X}Y_x}Z_{(x,y) }$ and $\coprod_{x \in X}\coprod_{y \in Y_x}Z_{(x,y)}$.

Now, for the monad unit, suppose we have an object $(X,\omega, (M_x)_{x \in X}) \in T^{'}A$, let us have a look at the following diagram:

\[\begin{tikzcd}
	{(X,\omega, (M_x)_{x \in X})} &&& {(X,\omega,(\{x\},*,M_x))} \\
	\\
	{(*,*,(X,\omega,(M_x)_{x \in X}))} &&& {(\coprod_{*}X,\omega,(M_x))= (\coprod_{x \in X}*,\omega,(M_x))=(X,\omega,(M_x))}
	\arrow[maps to, from=1-1, to=3-1]
	\arrow[maps to, from=3-1, to=3-4]
	\arrow[maps to, from=1-1, to=1-4]
	\arrow[maps to, from=1-4, to=3-4]
\end{tikzcd}\]

    Here the last equality in the two diagrams is, of course, an equality of the ordinal numbers which are order isomorphic to $\coprod_{*}X$,  $\coprod_{x \in X}*$ in the second diagram.

This shows that $T^{'}$ is a strict $2$-monad, and hence $T$ is a pseudo-monad, and their different categories of algebras are equivalent. Now, suppose that $A$ is a category, then we claim that the data of a normal colax algebra for the pseudo-monad defined above is the same as the data of an ultrastructure on $A$ as defined in \cite[definition 1.3.1]{lurie2018ultracategories}.

\section{From ultrastructures to colax algebras} \label{First direction}

\begin{theorem} \label{F} There exists a $2$-functorial process $\lu{F}$ from the bicategory of ultracategories, to the bicategory of $T$-normal colax algebras, and natural transformations, that assigns to an ultracategory a normal colax algebra for the pseudo-monad $T$, and such that this process makes the following triangle commutative, where $F_1, F_2$ are the forgetful functors: 
\[
\begin{tikzcd}[ampersand replacement=\&]
	{\mathrm{Ult}} \&\& {T\text{-}\mathrm{Alg}_{\mathrm{normal \ colax}}} \\
	\\
	\& {\mathrm{CAT}}
	\arrow["{\lu{F}}", from=1-1, to=1-3]
	\arrow["{F_1}"', from=1-1, to=3-2]
	\arrow["{F_2}", shift right=2, from=1-3, to=3-2]
\end{tikzcd}
\]
\end{theorem}
This means that $\lu{F}$ preserves the underlying categories, functors and natural transformations, and just assigns new structures to these.

In this section, we construct the colax algebra structure associated to the ultrastructure on a category $A$. Functoriality will be explained in \ref{Functoriality} and \ref{natural transformations}.

So suppose that $(A,(\int - d\omega)^A,\Delta^A,\epsilon^A)$  is an ultracategory (we will be often omitting the superscripts). We wish to construct the corresponding normal colax algebra structure for the pseudo-monad $T$ on $A$. We define the colax algebra functor $m^A$ on objects by sending $(X,\omega, (M_x)_{x \in X} ) $  to $\int_X M_x d \omega $.

Now on morphisms, suppose that we have $(X,f\omega,(M_x)_{x \in X})$ and $(Y,\iota_{Y^{'}}\omega,(N_y)_{y \in Y})$ and a morphism $ (f:Y^{'}\subseteq Y \rightarrow X ,(\chi_y)_{y \in Y^{'}})$ then the resulting map on morphisms is:

\begin{tikzcd}
\int_X M_x d f \omega \arrow[rr, "{\Delta_{\omega , f }}"] &  & \int_{Y^{'}} M_{f(y)} d  \omega  \arrow[rr, "\int \chi_y d\omega"] &  & \int_{Y^{'}} N_y d \omega \arrow[r, "\simeq"] & \int_Y N_y d \iota_{Y^{'}} \omega 
\end{tikzcd}

Here, the last isomorphism is the inverse of the ultraproduct diagonal map, which by axiom $(B)$ of \cite[definition 1.3.1]{lurie2018ultracategories} is invertible when the map (in this case, $\iota_{Y^{'}}$ is the inclusion of $Y^{'}$ in $Y$) is injective. Now we need to make sure our definition is independent of the representative of the class of the morphism, to do so suppose that we have $Y^{''} \subseteq Y^{'} \subseteq Y$, we are going to denote by $\iota_{Y^{''}}$ and $\iota_{Y^{'}}$ the respective inclusion maps in $Y$, and  $\iota_{Y^{''} \rightarrow Y^{'}}$ the inclusion of $Y^{''}$ in $Y^{'}$, and furthermore suppose that we have families of objects $(M_x)_{x \in X}$ and $(N_y)_{y \in Y}$, a  function $f$ from $Y^{'}$ to $X$ , and suppose that we have a family of morphisms $(\chi_y)_{y \in Y^{'}}$ from $M_{f(y)}$ to $N_y$. Now, let us have a look at the following diagram:

\adjustbox{max width =\textwidth}{
\begin{tikzcd}
	{\int_XM_xdf^{'}\omega=\int_XM_xdf\iota_{Y^{''}\rightarrow Y^{'}}\omega} &&& {\int_{Y^{''}}M_{f(y)} d\omega} && {\int_{Y^{''}}N_yd\omega} &&& {\int_{Y}N_y d\iota_{Y^{''}}\omega} \\
	\\
	&&& {\int_{Y^{'}}M_{f(y)} d \iota_{Y^{''}\rightarrow Y^{'}}\omega} && {\int_{Y^{'}}N_yd\iota_{Y^{''}\rightarrow Y^{'}}\omega}
	\arrow["{\Delta_{\omega,f \circ \iota_{Y^{''}\rightarrow Y'}}}", from=1-1, to=1-4]
	\arrow["{\Delta_{ \iota_{Y^{''}\rightarrow Y^{'}} \omega,f }}"', from=1-1, to=3-4]
	\arrow["{\int_{Y^{'}}\chi_yd\omega}", from=3-4, to=3-6]
	\arrow["{\int_{Y^{''}}\chi_yd\omega}", from=1-4, to=1-6]
	\arrow["{\Delta_{\omega,\iota_{Y^{''}\rightarrow Y^{'}}}}"', from=3-6, to=1-6]
	\arrow["{(\Delta_{\iota_{Y^{''}\rightarrow Y^{'}} \omega, \iota_{Y^{'}}})^{-1}}"', from=3-6, to=1-9]
	\arrow["{(\Delta_{\omega,\iota_{Y^{''}}})^{-1}}"{description}, from=1-6, to=1-9]
	\arrow["{\Delta_{\omega,\iota_{Y^{''}\rightarrow Y^{'}}}}"', from=3-4, to=1-4]
\end{tikzcd}
}
The left and the rightmost triangles commute \cite[proposition 1.3.5]{lurie2018ultracategories}, and the middle square is just naturality of the ultraproduct diagonal map. So the outermost diagram commutes, which results in the fact that the colax algebra functor is well-defined on morphisms.

Now we define the colax associator $a^{A}$ as follows: $ a^{A}_{(S,\omega,(T_s)_{s \in S}, (\gamma_s)_{s \in S},(M_{(s,t)})_{(s,t) \in \coprod_{s \in S} T_s})} =$

\[\begin{tikzcd}
	{\int_{\coprod_{s \in S}T_s}M_{(s,t)}d \int_S \iota_s \gamma_s d \omega } && {\int_S\int_{\coprod_{s \in S}T_s} M_{(s,t)} d \iota_s \gamma_s d\omega} && {\int_S\int_{T_s} M_{(s,t)} d \gamma_s d\omega}
	\arrow["{{\Delta_{\omega,\iota \gamma_{\bullet}}}}", from=1-1, to=1-3]
	\arrow["{{\int_S \Delta_{\gamma_s,\iota_s}}d\omega}", from=1-3, to=1-5]
\end{tikzcd}\]

Here $\iota_s$ is the inclusion map of $T_s$ inside $\coprod_{ s \in S} T_s$.

We will be omitting the subscript and the superscript of $a$ whenever it is convenient.

Let us show that $a$ is natural, to do so, let us explain the setting first:

Suppose that we have two objects of $T^2A$, $(S,\xi,(T_s, (M_{(s,t)}), \nu_s)_{s \in S})$ and $(Y,\omega^{'}, (X_y , N_{(y,x)}, \gamma^{'}_y)_{y \in Y})$, and a morphism between them. The data of such morphism is given by the following: 

\begin{itemize}
 
    \item a function  $f$ from $Y^{'}$ to $S$ such that $Y^{'}$ is a subset of $Y$ such that $Y^{'} \in \omega^{'}$ and such that if we call $\omega$ the restriction of $\omega'$ to $Y^{'}$  i.e. $\omega^{'} =i\omega$ (here $i$ is the inclusion of $Y^{'}$ in $Y$) , we have $\xi =f \omega$. 

    \item A family of functions $(g_y)_{y \in Y^{'}}$ defined on subsets $(X^{'}_y)_{y \in Y^{'}}$ such that for each of these subsets $X^{'}_y \subseteq X_y$ satisfies $X^{'}_y \in \gamma^{'}_y$, and if we call $\gamma_y$ the restriction of  $\gamma^{'}_y$ to $X^{'}_y$ i.e. $\gamma^{'}_y =\iota_y \gamma_y$ (here $\iota_y$ is the inclusion map of $X^{'}_y$ inside $X_y$),  then we have $\nu_{f(y)} = g_y \gamma_y$.

    \item A family of morphisms in the category $A$ indexed by $\coprod_{y \in Y'}X'_y$, which we will call $(\psi_{(y,x)})$, such that $\psi_{(y,x)}: M_{(f(y),g_y(x))} \rightarrow N_{(y,x)}$.
    
\end{itemize}

We are going to call $\bar{f}$ the map from $\coprod_{y \in Y^{'}} X^{'}_y$ to $\coprod_{s \in S} T_s$ which sends $(y,x)$ to $(f(y), g_y(x))$, and we are going to call $j$ the inclusion function of $\coprod_{y \in Y^{'}}X^{'}_y$ in $\coprod_{y \in Y}X_y$. Naturality of $a$ translates to the commutativity of the outermost square in \hyperref[diag1]{Diagram 1} of \hyperref[appendix2]{Appendix B}.

Let us have a look at this diagram: \\
\noindent Squares $2$,$4$,$5$, and $7$ are naturality squares.\\
\noindent Squares $3$ and $8$ commute by \cite[proposition 1.3.5]{lurie2018ultracategories}.\\
\noindent So the only thing remaining is showing that squares $1$ and $6$ are commutative in \hyperref[diag1]{Diagram 1} of \hyperref[appendix2]{Appendix B}.

Starting with square $1$, we use the following diagram, we want to show that the outermost square is commutative:
\[
\adjustbox{scale =0.77}{
\begin{tikzcd}[ampersand replacement=\&]
	{\int_{\coprod_{s \in S}T_s}M_{(s,t)}d \int_S \iota_s \nu_s d f\omega } \&\&\& {\int_S\int_{\coprod_{s \in S}T_s} M_{(s,t)} d \iota_s \nu_s df\omega} \\
	\\
	\\
	\&\&\& {\int_{Y^{'}} \int_S \int_{\coprod_{s \in S}T_s}M_{(s,t)}d \iota_{s}\nu_s d \delta_{f(y)}d\omega} \\
	\\
	\\
	{\int_{\coprod_{y \in Y^{'}}X^{'}_y}M_{(f(y), g_y(x))}d\int_{Y^{'}}\iota_{y} \gamma_{y}d\omega} \&\&\& {\int_{Y^{'}} \int_{\coprod_{s \in S}T_s}M_{(s,t)}d \iota_{f(y)}g_{y}\gamma_{y}d\omega} \\
	\\
	\\
	\\
	\\
	\&\&\& {\int_{Y^{'}}\int_{\coprod_{y \in Y^{'}}X^{'}_y}M_{(f(y), g_y(x))}d\iota_{y} \gamma_{y}d\omega}
	\arrow["{{\Delta_{f\omega,\iota \nu_{\bullet}}}}", from=1-1, to=1-4]
	\arrow["{\Delta_{\int_{Y^{'}} \iota_y \gamma_y d \omega , \bar{f}}}"', from=1-1, to=7-1]
	\arrow[""{name=0, anchor=center, inner sep=0}, "{\Delta_{\omega,\iota \gamma_{\bullet}}}"{description}, from=1-1, to=7-4]
	\arrow[""{name=1, anchor=center, inner sep=0}, "{ \Delta_{\omega,\delta_{f(y)}}}", from=1-4, to=4-4]
	\arrow["{\int_{Y^{'}}\epsilon_{S,f(y)} d\omega =\int_{Y'} \Delta^{-1}_{\delta_{f(y)},\iota\nu_{\bullet}}d\omega}"{description}, from=4-4, to=7-4]
	\arrow[""{name=2, anchor=center, inner sep=0}, "{{\Delta_{\omega,\iota \gamma_{\bullet}}}}"', from=7-1, to=12-4]
	\arrow["{\int_{Y^{'}}\Delta_{ \iota_{y}\gamma_{y},\bar{f}} d\omega}"{description}, from=7-4, to=12-4]
	\arrow["1"{description}, draw=none, from=0, to=1]
	\arrow["2"{description}, draw=none, from=0, to=2]
\end{tikzcd}
}
\]

Square $1$ commutes by Axiom $C$ of \cite[definition 1.3.1]{lurie2018ultracategories} of  ultracategory axioms, while square $2$ is just a naturality square.

Now, for square $6$ of  \hyperref[diag1]{Diagram 1} of \hyperref[appendix2]{Appendix B} , its commutativity  corresponds to the commutativity of the outermost square of the following diagram:

\adjustbox{scale =0.65}{
\begin{tikzcd}[ampersand replacement=\&]
	{\int_{\coprod_{y \in Y}X_y}N_{(y,x)} d \int_{Y} \iota_y \gamma_{y}^{'}d\omega^{'}} \&\&\& {\int_{Y}\int_{\coprod_{y \in Y }X_y}N_{(y,x)}d\iota_{y} \gamma^{'}_{y}d\omega^{'}} \&\&\&\& {\int_Y \int_{\coprod_{y \in Y}X_y} N_{(y,x)}d\iota_y\gamma_y d\omega^{'}} \\
	\\
	\& {\int_{Y^{'}}\int_{\coprod_{y \in Y}X_y}N_{(y,x)} d \iota_y\gamma^{'}_yd\omega} \&\& {\int_{Y^{'}}\int_{Y}\int_{\coprod_{y\in Y}X_y}N_{(y,x)}d\iota_y\gamma^{'}_yd \delta_{y}d\omega} \&\&\&\& {\int_{Y^{'}}\int_{Y}\int_{\coprod_{y\in Y^{'}}X^{'}_y}N_{(y,x)}d\iota_y\gamma_yd \delta_{y}d\omega} \\
	\\
	\\
	\\
	\\
	{\int_{\coprod_{y \in Y^{'}}X^{'}_y}N_{(y,x)}d\int_{Y^{'}}\iota_{y} \gamma_{y}d\omega} \&\&\&\&\&\&\& {\int_{Y^{'}}\int_{\coprod_{y \in Y^{'}}X^{'}_y}N_{(y,x)}d\iota_{y} \gamma_{y}d\omega}
	\arrow["{{{{\Delta_{\omega^{'},\iota \gamma^{'}_{\bullet}}}}}}"', from=1-1, to=1-4]
	\arrow["{{{\Delta_{\omega,\iota\gamma^{'}_{\bullet}}}}}", from=1-1, to=3-2]
	\arrow["2"{description}, draw=none, from=1-1, to=3-4]
	\arrow["{{{\Delta_{\int_{Y'}\iota_y\gamma_yd\omega,j}}}}"', from=1-1, to=8-1]
	\arrow["{{{\int_{Y} \Delta_{\iota_y\gamma_y,\iota_y}d \omega^{'}}}}", from=1-4, to=1-8]
	\arrow["{{{\Delta_{\omega,\delta_{(\bullet)}}}}}"', from=1-4, to=3-4]
	\arrow["3", draw=none, from=1-4, to=3-8]
	\arrow["{{{\Delta_{\omega, \delta_{\bullet}}}}}", from=1-8, to=3-8]
	\arrow["{{{\int_{Y^{'}}\Delta_{\iota_y\gamma_y,j}d\omega}}}"', from=3-2, to=8-8]
	\arrow[""{name=0, anchor=center, inner sep=0}, "{{{\int_{Y^{'}}\epsilon_{Y,y}d\omega}}}"', from=3-4, to=3-2]
	\arrow["{{{\int_{Y^{'}}\int_{Y} \Delta_{\iota_y\gamma_y,j}d \delta_y d\omega}}}"{description}, from=3-4, to=3-8]
	\arrow[""{name=1, anchor=center, inner sep=0}, "{{{\int_{Y^{'}}\epsilon_{Y,y} d\omega}}}", from=3-8, to=8-8]
	\arrow[""{name=2, anchor=center, inner sep=0}, "{{{{\Delta_{\omega,\iota \gamma_{\bullet}}}}}}", from=8-1, to=8-8]
	\arrow["1", draw=none, from=3-2, to=2]
	\arrow["4"', draw=none, from=0, to=1]
\end{tikzcd}
}

Square $2$ is commutative by \cite[Axiom C of definition 1.3.1]{lurie2018ultracategories}, while squares $3$ and $4$ are naturality squares, it remains to check the commutativity of square $1$ of the diagram above.

To do so, let us use the following diagram, as usual we want to show that the outermost diagram is commutative:
\[
\adjustbox{scale ={0.55}{0.7}}{

\begin{tikzcd}
	{\int_{\coprod_{y \in Y}X_y}N_{(y,x)} d \int_{Y} \iota_y \gamma_{y}^{'}d\omega^{'}=\int_{\coprod_{ y \in Y}X_y}N_{(y,x)}d \int_{\coprod_{y \in Y^{'}}X^{'}_y}\delta_{(y,x)}d\int_{Y^{'}} {}\iota_{y}\gamma_{y}d\omega=\int_{\coprod_{y \in Y}X_y}N_{(y,x)} d \int_{Y^{'}} \iota_y \gamma_{y}^{'}d\omega} &&&& {\int_{Y^{'}}\int_{\coprod_{y \in Y}X_y}N_{(y,x)} d \iota_y\gamma^{'}_yd\omega=\int_{Y^{'}}  \int_{\coprod_{y \in Y}X_y} N_{(y,x)}d \int_{\coprod_{y \in Y^{'}}X^{'}_y} \delta_{(y,x)}d\iota_y \gamma_y d\omega} \\
	& {} \\
	{\int_{\coprod_{y \in Y^{'}}X^{'}_y} \int_{\coprod_{y \in Y}X_y} N_{(x,y)}d\delta_{(y,x)}d\int_{Y^{'}}\iota_y \gamma_{y}d\omega} &&&& {\int_{Y^{'}} \int_{\coprod_{y \in Y^{'}}X^{'}_y}\int_{\coprod_{ y \in Y}X_y}N_{(y,x)} d\delta_{(y,x)}d\iota_{y}\gamma_{y}d\omega} \\
	\\
	{\int_{\coprod_{y \in Y^{'}}X^{'}_y}N_{(y,x)}d\int_{Y^{'}}\iota_{y} \gamma_{y}d\omega} &&&& {\int_{Y^{'}}\int_{\coprod_{y \in Y^{'} }X^{'}_y}N_{(y,x)}d\iota_{y} \gamma_{y}d\omega}
	\arrow[""{name=0, anchor=center, inner sep=0}, "{\Delta_{\omega,\iota\gamma^{'}_{\bullet}}=\Delta_{\omega,\int_{\coprod_{y \in Y^{'}}X^{'}_y}\delta_{(y,x)}d\iota\gamma_{\bullet}}}", from=1-1, to=1-5]
	\arrow["{\Delta_{\int_Y\iota_y\gamma_yd\omega,\delta_{\bullet}}}"', from=1-1, to=3-1]
	\arrow["{\int_{Y^{'}}\Delta_{\iota_y\gamma_y,j}d\omega}", from=1-5, to=3-5]
	\arrow[""{name=1, anchor=center, inner sep=0}, "{\Delta_{\omega,\iota\gamma_{\bullet}}}"', from=3-1, to=3-5]
	\arrow["{\int_{\coprod_{y \in Y^{'}}X^{'}_y} \epsilon_{\coprod_{y\in Y}X_y,(y,x)} d\int_{Y^{'}} \iota_y \gamma_y d\omega}"', from=3-1, to=5-1]
	\arrow["{\int_{Y^{'}} \int_{\coprod_{y \in Y^{'}}X^{'}_y} \epsilon_{\coprod_{y\in Y}X_y,(y,x)}d \iota_{y} \gamma_y d\omega}", from=3-5, to=5-5]
	\arrow[""{name=2, anchor=center, inner sep=0}, "{{\Delta_{\omega,\iota \gamma_{\bullet}}}}"', from=5-1, to=5-5]
	\arrow["1", draw=none, from=0, to=1]
	\arrow["2", draw=none, from=1, to=2]
\end{tikzcd}
}
\]

Square $1$ commutes by axiom $C$ of the ultracategory axioms of \cite[definition 1.3.1]{lurie2018ultracategories}, while square $2$ commutes by naturality of ${\Delta_{\omega,\iota \gamma_{\bullet}}}$.

A particular case of this definition is when we have the same set $T$ repeated $|S|$ times, and a family of objects $(M_t)_{t \in T}$ indexed by $T$, in the following we define the family $(M^{'}_{(s,t)})_{(s,t) \in S \times T}$  as the family such that $M^{'}_{(s,t)}=M_t$.
\[
\begin{tikzcd}
{\int_{S \times T}M^{'}_{(s,t)}d \int_S \iota_s \gamma_s d \omega } \arrow[rr, "{\Delta_{\omega,\iota \gamma_{\bullet}}}"] &  & {\int_S\int_{S \times T} M^{'}_{(s,t)} d \iota_s \gamma_s d\omega} \arrow[rr, "{\int_S \Delta_{\gamma_s,\iota_s}d\omega}"] &  & \int_S\int_T M_t d \gamma_s d\omega
\end{tikzcd}\]

Finally, we define the colax unitor $i^{A}_{M}$ by:

 \[\begin{tikzcd}
\int_* M d* \arrow[rr, "\epsilon_{*,*}"] &  & M
\end{tikzcd}\]
which is hence an isomorphism. Naturality in $M$ is immediate. So formally, we define \newline $\lu{F}((A,(\int-d\omega)^A,\Delta^A,\epsilon^A))=(A,m^A,a^{A},i^{A})$.

Now, we want to show that this is really a colax algebra for the pseudo-monad $T$, towards that, we are going to show that this definition satisfies the colax algebra axioms \cite[section 3]{penon2009representable} when restricted to ordinals, so we define $T^{'}$ to be the restriction of $T$ to sets which are ordinals (so $T' A$ is the full subcategory of $TA$ whose objects are triples $(X,\omega, (M_x)_{x \in X})$ where $X$ is an ordinal), and hence by definition the following diagram commutes:
\[\begin{tikzcd}
	{T^{'}A} && A \\
	TA
	\arrow["{m^{'}}", from=1-1, to=1-3]
	\arrow["\simeq"', from=1-1, to=2-1]
	\arrow["m"', from=2-1, to=1-3]
\end{tikzcd}\]

We define ${a^{'}}^A$ and ${i^{'}}^A$ to be the corresponding restrictions of $a^A$ and $i^A$.

Now we turn to checking that our definition satisfies the axioms of colax algebra; Starting with the first equality which is the following:

\begin{tikzcd}[ampersand replacement=\&]
	T'A \&\& {{T'}^2A} \&\& T'A \&\&\& {=} \& T'A \&\&\& A \\
	\&\& A \&\& T'A \&\& A
	\arrow["{\eta_{T'A}}", from=1-1, to=1-3]
	\arrow["m'"', from=1-1, to=2-3]
	\arrow["{\mu_A}", from=1-3, to=1-5]
	\arrow["T'm'"{description}, from=1-3, to=2-5]
	\arrow["a'"', Rightarrow, nfold, from=1-5, to=2-5]
	\arrow["m'", from=1-5, to=2-7]
	\arrow["m'"', from=1-9, to=1-12]
	\arrow["{\eta_A}"', from=2-3, to=2-5]
	\arrow[""{name=0, anchor=center, inner sep=0}, "{\la{Id}}"', curve={height=30pt}, from=2-3, to=2-7]
	\arrow["m'", from=2-5, to=2-7]
	\arrow["i"{description}, between={0}{0.8}, Rightarrow, nfold, from=2-5, to=0]
\end{tikzcd}

This corresponds to the commutativity of the following diagram :
\[
\begin{tikzcd}
\int_I M_{i} d \omega\arrow[dd, "{a^{'}_{(*,*,(I),  \omega, (M_{i})_{ i \in I})}=\Delta_{*,\omega}}"'] \arrow[rrdd, "id"] &  &                 \\
                                                                                              &  &                 \\
\int_*\int_I M_i d \omega d * \arrow[rr, "\epsilon_{*,*}"']                                            &  & \int_I M_i d\omega
\end{tikzcd}\]

Notice that $a^{'}_{(*,*,(I),  \omega, (M_{i})_{ i \in I})}=\Delta_{*,\omega}$ requires us to use the equality of the ordinals $\coprod_{*}I$ (which means the unique ordinal order isomorphic to $\coprod_{*}I$ where we equip the last with the lexicographic order) and $I$. Here $*$ is the one-point set. This diagram commutes by axiom $A$ of Lurie's ultracategory axioms.

Now we show the second equality which is the following:

\begin{tikzcd}[ampersand replacement=\&]
	T'A \&\& {{T'}^2A} \&\& T'A \&\& A \& {=} \& T'A \&\&\& A \\
	\&\&\&\& T'A
	\arrow["{T'\eta_A}", from=1-1, to=1-3]
	\arrow[""{name=0, anchor=center, inner sep=0}, "{\la{Id}}"', curve={height=18pt}, from=1-1, to=2-5]
	\arrow["{\mu_A}", from=1-3, to=1-5]
	\arrow["T'm'"', from=1-3, to=2-5]
	\arrow["m'", from=1-5, to=1-7]
	\arrow["a'"', Rightarrow, nfold, from=1-5, to=2-5]
	\arrow["m'", from=1-9, to=1-12]
	\arrow["m'"', from=2-5, to=1-7]
	\arrow["T'i"', between={0}{0.85}, Rightarrow, nfold, from=1-3, to=0]
\end{tikzcd}

this corresponds to the commutativity of the following diagram:
\[
\begin{tikzcd}
	{\int_I M_{i} d \omega} \\
	\\
	{\int_I\int_{I} M_i d \delta_i d\omega} && {\int_I M_{i} d \omega}
	\arrow["{{\Delta_{\omega,\delta_{\bullet}}=a^{'}_{(I,\omega,(\{i\})_{i \in I},(\delta_i)_{i \in \{i\}},(M_i)_{i \in \{i\}})}}}"', from=1-1, to=3-1]
	\arrow["id", from=1-1, to=3-3]
	\arrow["{{\int \epsilon_{I,i} d \omega}}", from=3-1, to=3-3]
\end{tikzcd}\]

But this is actually one of the results by Lurie \cite[corollary 1.3.6]{lurie2018ultracategories}. Here we use the equality of the two ordinal numbers  $\coprod_{i\in I}*$ (which again means the unique ordinal order isomorphic to $\coprod_{i \in I}*$ where we equip the last with the lexicographic order) with $I$. 

Finally, we check the third axiom:

\[\begin{tikzcd}
	{m^{'} \circ \mu_A \circ \mu_{T^{'}A}} &&& {m^{'} \circ \mu_A \circ T^{'} \mu_A} &&& {m^{'} \circ T^{'} m^{'} \circ T^{'} \mu_A} \\
	\\
	{m^{'} \circ T^{'} m^{'} \circ \mu_{T^{'}A}} &&& {m^{'} \circ \mu_{A} \circ {T^{'}}^2 m^{'}} &&& {m^{'} \circ T^{'} m^{'} \circ {T^{'}}^2 m^{'}}
	\arrow["{{{{a^{'} \circ \mu_{ T'A}}}}}"', from=1-1, to=3-1]
	\arrow["{{{{\mathrm{Id}}}}}"', from=1-4, to=1-1]
	\arrow["{{{{a^{'} \circ T' \mu_A}}}}", from=1-4, to=1-7]
	\arrow["{{{{m^{'} \circ T^{'}a^{'}}}}}", from=1-7, to=3-7]
	\arrow["{{{{\mathrm{Id}}}}}"', from=3-1, to=3-4]
	\arrow["{{{{a^{'} \circ {{T}^{'}}^2 }}}}"', from=3-4, to=3-7]
\end{tikzcd}\]

This corresponds to the commutativity of \hyperref[diag2]{Diagram 2} of \hyperref[appendix2]{Appendix B}: Here $\rho$ in \hyperref[diag2]{Diagram 2} is equal to $\int_R   \iota_r (\int_{S_r} \iota_s \nu_{(r,s)} d \omega_r) d\lambda=\int_R    (\int_{\coprod_{R} S_r}\iota_r\iota_s \nu_{(r,s)} d \iota_r \omega_r) d\lambda=\int_{\coprod_{R}S_r}\iota_r \iota_s \nu_{(r,s)}d \int_R  \iota_r \omega_r d\lambda$.

Now to show this commutativity, we use \hyperref[diag3]{Diagram 3} of \hyperref[appendix2]{Appendix B}:

\noindent Square $1$ of \hyperref[diag3]{Diagram 3} commutes by axiom $(C)$ of Lurie's axioms.\\
Square $3$ commutes by naturality of $\Delta_{\lambda,\iota\omega_{\bullet}}$.\\
Square $4$ commutes by naturality of the ultraproduct diagonal maps $\Delta_{\omega_r,\iota_r}$.\\
Triangle $6$ commutes by proposition $1.3.5$ in Lurie's paper.\\
So, the only two diagrams that are a little bit problematic are triangle $2$ and square $5$.\\
Square $5$ commutes by axiom $(C)$ of Lurie's axioms and by naturality of $\epsilon_{S \times T,\delta_{\iota_r(t)}}$ (we leave the technical details).\\
Finally, for triangle $2$ we need  proposition $1.3.5$ and axiom $(C)$. The proof technique is the same as in section \nameref{3}, but we are in a much simpler case.\\
The diagram in the proof should look like:

\adjustbox{max width = \textwidth}{

\begin{tikzcd}
	{ \int_{\coprod_{R}\coprod_{S_r}T_{(r,s)}} Md \int_{\coprod_{R}S_r} \iota_{(r,s)}\nu_{(r,s)}d \int_{S_r} \delta_{(r,s)} d\omega_r  } &&&& {\int_{\coprod_{R}S_r}\int_{\coprod_{R}\coprod_{S_r}T_{(r,s)}}M d\iota_{(r,s)}\nu_{(r,s)}d\iota_r\omega_r} \\
	&&&& {} \\
	\\
	{\int_{S_r} \int_{\coprod_{R}\coprod_{S_r}T_{(r,s)}} Md \iota_r \iota_s\nu_{(r,s)}d\omega_r} &&&& {\int_{S_r} \int_{\coprod_{R}S_r} \int_{\coprod_{R}\coprod_{S_r}T_{(r,s)}}M d\iota_{(r,s)}\nu_{(r,s)}d\delta_{(r,s)}d\omega_r} \\
	\\
	{} && {} & {} & {\int_{S_r}\int_{\coprod_{R}\coprod_{S_r}T_{(r,s)}} Md  \iota_{s}\nu_{(r,s)}d\omega_r}
	\arrow["{\Delta_{\iota_r\omega_r,\iota_{\bullet}\nu_{\bullet}}}", from=1-1, to=1-5]
	\arrow["{\Delta_{\omega_r,\iota_{\bullet}}}", from=1-1, to=4-1]
	\arrow["{\Delta_{\omega_r,\delta_{(r,\bullet)}}}", from=1-5, to=4-5]
	\arrow["{\int_{S_r}\Delta_{\delta_{(r,s)},\iota_{\bullet}}d\omega_r}", from=4-1, to=4-5]
	\arrow["{\la{Id}}", from=4-1, to=6-5]
	\arrow["{ \int_{S_r}\epsilon_{(r,s),\coprod_{R}S_r} d\omega_r}", from=4-5, to=6-5]
\end{tikzcd}
}
The upper diagram commutes by axiom $(C)$ of \cite[definition 1.3.1]{lurie2018ultracategories}, while the lower diagram commutes using \cite[corollary 1.3.6]{lurie2018ultracategories}, finally notice that the  composition $ \int_{S_r}\epsilon_{(r,s),\coprod_{R}S_r} d\omega_r \circ \Delta_{\omega_r,\delta_{(r,\bullet)}} $   is exactly  the definition of $\Delta_{\omega_r,\iota_r}$. Here, $M$ is just an abbreviation for $M_{(r,s)}$.

\section{From colax algebras to ultrastructures} \label{Second direction}

\begin{theorem}
There exists a $2$-functorial process $\lu{G}$ that sends a category with a normal colax algebra structure for the pseudo-monad $T$, to the same category with an ultrastructure, and such that this $2$-functor makes the following diagram commutative: 

\[
\begin{tikzcd}[ampersand replacement=\&]
	{T\text{-}\mathrm{Alg}_{\mathrm{normal \ colax}}} \&\& {\mathrm{Ult}} \\
	\\
	\& {\mathrm{CAT}}
	\arrow["{\lu{G}}", from=1-1, to=1-3]
	\arrow["{F_2}"', from=1-1, to=3-2]
	\arrow["{F_1}", shift right=2, from=1-3, to=3-2]
\end{tikzcd}
\]
\label{G}
\end{theorem}

This means that $\lu{G}$ does not only preserve the category, but also does not change functors and natural transformations, it just assigns new structure to those.

Again, in this section, we construct the ultrastructure associated to the $T$-normal colax algebra structure on a category $A$. Functoriality will be explained in \ref{Functoriality} and \ref{natural transformations}.

Suppose that a category $A$ has the structure of $T$-normal colax algebra, we denote by $m^A$ the algebra morphism from $TA$ to $A$, which is again, as indicated by our definition, a functor \\ $m^{A}: TA \rightarrow A$, together with a normal colax algebra structure ${m^{'}}^{A}: T^{'}A \rightarrow A $ (with colax associator ${a^{'}}^A$ and unitor ${i^{'}}^A$) such that the following diagram commutes:
\[\begin{tikzcd}
	{T^{'}A} && A \\
	TA
	\arrow["{{m^{'}}^{A}}", from=1-1, to=1-3]
	\arrow["\simeq"', from=1-1, to=2-1]
	\arrow["m^{A}"', from=2-1, to=1-3]
\end{tikzcd}\]

Here, the equivalence above is the full subcategory inclusion of $T^'A$ inside $TA$. As usual we will be omitting superscripts whenever it's convenient.

Before continuing, we must note that ${a^'}^A$ can be extended  from a map between $m^{'} \circ \mu_A$ and  $m^{'} \circ T'm^{'}$ , to a map $a^A$ from $m \circ \mu_A$ to $m \circ Tm$. Similarly, we extend ${i^{'}}^A$  from a map from $m^{'} \circ \eta_{A}$ to $\la{Id}$ to a map $i^{A}$ from $m \circ \eta_{A}$ to $\la{Id}$, in the following way:

\paragraph{Colax associator and unitor of $T$ algebra structure}
\label{a}

Let us define a $2$-morphism $a$ from  $m \circ \mu_A$ to $m \circ Tm$  which extends $a^{'}$, suppose that we have a family of objects  $(I_s)_{s \in S}$ as well as an ultrafilter $\omega$ on $S$, then if we choose a well-ordering on each of the $I_s$ as well as a well-ordering on $S$, and a family of ultrafilters on each $I_s$, which we call $(\lambda_s)_{s \in S}$; let us denote by $\phi(X)$ the corresponding ordinal for a set $X$, then we can form $a_{S,\omega,(I_s)_{s \in S}, (\lambda_s)_{s \in S},(M_{(s,i)})_{(s,i) \in \coprod_{s \in S} I_s}}$ as follows:

\[\begin{tikzcd}
	{\int_{\coprod_{S}I_s} M_{(s,i)}d \int_S \iota_s\lambda_s d\omega} &&&& {} & {\int_S \int_{I_s}M_{(s,i)}d\iota_s \lambda_s d\omega} \\
	\\
	{\int_{\coprod_{\phi(S)}\phi(I_s)}M_{(s,i)}d e_1\int_S \iota_s\lambda_s d\omega} && {} && {} & {\int_{\phi(S)} \int_{\phi(I_s)}M_{(s,i)}d\iota_s e_s\lambda_s de_2\omega}
	\arrow["{a_{S,\omega,(I_s)_{s \in S}, (\lambda_s)_{s \in S},(M_{(s,i)})_{(s,i) \in \coprod_{s \in S} I_s}}}", from=1-1, to=1-6]
	\arrow["{m(e_1)}"', from=1-1, to=3-1]
	\arrow["{a^{'}_{(\phi(S),\omega,(\phi(I_s))_{s \in \phi(S)}, (\lambda_s)_{s \in \phi(S)},(M_{(s,i)})_{(s,i) \in \coprod_{s \in \phi(S)} \phi(I_s)})}}", from=3-1, to=3-6]
	\arrow["{(m \circ Tm)(e_2, (e_s)_{s \in S})}"', from=3-6, to=1-6]
\end{tikzcd}\]

Here $e_1, e_2 , (e_s)_{s \in S}$ are the isomorphisms between $\coprod_{\phi(S)}\phi(I_s)$ and $\coprod_{S}I_s$, $S$ and $\phi(S)$, $I_s$ and $\phi(I_s)$ respectively. Showing that our construction is natural is fairly easy, we can also show that it is independent of the ordinal choice, although we are not going to be needing this. So we were able to extend $a^{'}$ to a  $2$-morphism $a$ from  $m \circ \mu_{A}$ to $m \circ Tm$.

We may also define $i_{M} =i^{'}_{M}$ for any object $M \in A$ (assume for simplicity that we take a canonical choice of one-point set to be the one-point ordinal, otherwise we may construct $i$ similarly). And hence, by virtue of this construction, we may consider $(A, m^A, a^A, i^A)$ to be the data of $T$-colax algebra (if we had used the actual definition of a $T$-colax algebras, this would be really the data of colax algebras, but we have used an equivalent definition where we only restricted our attention to ordinals, but we actually need the data in its full generality for our next construction, so we had to recover it first).

We define the ultracategory structure associated to $(A, m^A, a^A, i^A)$  as follows: Let $*$ simultaneously denote the one-point set and the unique ultrafilter on that set, as usual. A morphism in $TA$ is a pair $(f,(b_y))$ as we have stated  before, now when the family $(b_y)$ is just identity maps we are going to omit writing it. The ultracategory functor is defined as follows :\\
For every set $I$ and every ultrafilter $\omega$ on $I$, we define the functors $\kappa_{I,\omega}$ from $A^I$ to $TA$ that send $(M_i)_{i \in I}$ to $(I,\omega,(M_i)_{i \in I})$ and on morphisms it sends $(b_i)_{i \in I}$ to $(\la{Id},(b_i)_{i \in I})$. Then the ultraproduct functor corresponding to the ultrafilter $\omega$ is defined as $m \circ  \kappa_{I,\omega}$, and as usual we are going to denote it by $(\int_{I} \bullet d\omega)^A.$

Suppose that $f$ is a map of sets between $X$ and $Y$ and let $\omega$ be an ultrafilter on $X$, and suppose that we have a family of objects $(M_y)_{y \in Y}$, then we can define a morphism $(f,(id_{M_{y}}))$ between $(Y,f\omega,(M_{y}))$ and  $(X, \omega,(M_{f(x)})_{x \in X})$ in the category $TA$, then we are going to write $m_{f}$ instead of writing $m_{(f,(id_{M_{y}}))}$  for the corresponding morphism from $\int_Y M_y df\omega$ to $\int_{X}M_{f(x)}d\omega$.

Now, for the categorical Fubini transform, we define it as follows:

\begin{tikzcd}
\int_TM_td(\int_S \gamma_s d \omega) \arrow[rr, "m_{\pi_T}"] &  & {\int_{S \times T} M^{'}_{(s,t)} d \int_S \iota_s \gamma_s d \omega} \arrow[rr, "{a_{S,\omega,(T^{'}_s),(\gamma_s),(M^{'}_{(s,t)})}}"] &  & \int_S \int_T M_t d \gamma_s d \omega
\end{tikzcd}

Here, for every $s \in S$, $T'_s$ is defined as $T^{'}_s =T$, and $M^{'}_{(s,t)}$ is defined as $M^{'}_{(s,t)}= M_t$. So we have $S\times T =\coprod_{s \in S} T^{'}_s$, and of course we define $\iota_s$ as the inclusion map of a copy of $T$ (namely $T'_s$) in $S \times T$.

\noindent Naturality in the family $(M_t)$ is obvious. Here  $a$ denotes the colax algebra associator of $T$ as defined above.

Now, we define the natural isomorphism $\epsilon_{S,s_0}$ as the following composition :

\begin{tikzcd}
	{\int_S M_sd\delta_{s_0}} &&& {\int_*M_{s_0} d*} &&& M_{s_0}
	\arrow["{m_{* \mapsto s_0}}", from=1-1, to=1-4]
	\arrow["{i_{M_{s_0}}}", from=1-4, to=1-7]
\end{tikzcd}

\noindent This map is an isomorphism, since it's a composition of two isomorphisms. Naturality in the family $(M_s)$ is again clear.

So we define $\lu{G}(A,m^A,a^A,i^A)= (A,(\int- d\omega)^A,\Delta^A,\epsilon^A)$.

Now we need to check that our definition satisfies ultracategory axioms of \cite[Definition 1.3.1]{lurie2018ultracategories}.

\subsection*{Axiom A} We want to show that our data  satisfies the first of Lurie's axiom, which means that the following composition is the identity: 
\[\begin{tikzcd}
	{\int_T M_t d(\int_S \gamma_sd \delta_{s_0})} &&& {\int_S \int_T M_t d \gamma_s d \delta_{s_0}} &&& {\int_T M_t d \gamma_{s_0}} \\
	{} \\
	{}
	\arrow["{\Delta_{\delta_{s_0},\gamma_{\bullet}}}"', from=1-1, to=1-4]
	\arrow["{\epsilon_{S,s_0}}"', from=1-4, to=1-7]
\end{tikzcd}\]


To do so, we use the following diagram:

\[\begin{tikzcd}
	{\int_T M_t d(\int_S \gamma_s d \delta_{s_0})} &&&& {\int_{T} M_t d( \int_* \gamma_{s_0}d*)} \\
	\\
	\\
	{\int_{S \times T} M_t d( \int_S \iota_s\gamma_s d\delta_{s_0})} &&&& {\int_{T} M_t d( \int_* \gamma_{s_0}d*)} \\
	\\
	&&&&& {} \\
	{\int_S \int_T M_t d \gamma_s d \delta_{s_0}} &&&& {\int_*\int_T M_t d \gamma_{s_0} d*} \\
	\\
	\\
	&&&& {\int_T M_t d \gamma_{s_0}}
	\arrow["{\la{Id}}", from=1-1, to=1-5]
	\arrow["{m_{\pi} }"', from=1-1, to=4-1]
	\arrow["{\la{Id}}", from=1-5, to=4-5]
	\arrow["{m_{t \mapsto (s_0,t)}}"', from=4-1, to=4-5]
	\arrow["{a_{(S,\delta_{s_0},(T^{'}_s),\gamma_s,(M^{'}_{(s,t)}))}}"', from=4-1, to=7-1]
	\arrow["{a_{(*,*,T,\gamma_{s_0},(M_t))}}"', from=4-5, to=7-5]
	\arrow["{m_{*\mapsto s_0}}"', from=7-1, to=7-5]
	\arrow["{\epsilon_{S,s_0}}"', from=7-1, to=10-5]
	\arrow["{i_{\int_TM_td\gamma_{s_0}}}"', from=7-5, to=10-5]
\end{tikzcd}\]

We wish to show that the outermost diagram commutes: the upper square commutes by functoriality of $m$, the middle square is naturality of $a$, and the lower triangle is the definition of $\epsilon_{S,s_0}$.

\subsection*{Axiom B}
Now we show that our construction satisfies axiom $B$ of Lurie's axioms, but this automatically follows from our definition of the monad, since suppose that we have an injective map $f$ from $X$ to $Y$ and some ultrafilter $\omega$ on $X$, then the map from  $(Y,f\omega ,(M_y))$ to $(X,\omega,(M_{f(x)}))$ is invertible, so its image by the functor $m$ is also invertible.

\subsection*{Axiom C} \label{axiom 3}

Now, to show that axiom $C$ of \cite[Definition 1.3.1]{lurie2018ultracategories} is satisfied by our definition: We are working in the same setting as axiom $C$ of Lurie's axioms, and we want to show that the diagram \hyperref[diag4]{Diagram 4} of \hyperref[appendix2]{Appendix B} is commutative, where $\rho$  is defined by $\int_R \int_S \nu_{s} d\omega_r d\lambda =\int_S \nu_s d (\int_R  \omega_r  d\lambda)$.

We are going to do the proof in the case where all sets are ordinals, and we leave the more general proof for arbitrary sets to  \hyperref[A]{Appendix A}. In this specific case, the maps $m$ and $m^{'}$ agree, also this holds for the maps $a$ and $a^{'}$. Now let us have a look at \hyperref[diag5]{Diagram 5} of \hyperref[appendix2]{Appendix B}, we have that:

\noindent Square $1$ commutes by functoriality of $m^{'}$.\\
Squares $2,3$ commute by naturality of the associator $a^{'}$.\\
Square $4$ commutes by the second diagram of the definition of colax algebras as defined in  \cite[section 3]{penon2009representable}, and hence the outermost diagram commutes, which is exactly what we wanted to show.

\section{The two processes are the inverses of each other} \label{Showing bijection}

Our goal now is to show that the two processes ($\lu{F}$ and $\lu{G}$) of \ref{F} and \ref{G} are actually the inverses of each other:

The first thing we should notice is that the two processes we just showed do not affect the underlying categories as well as the ultraproduct and the algebra functors, so for an ultracategory $A$ we have the equality of $(\int - d\omega)^A$ and $(\int - d\omega)^{\lu{G}(\lu{F}(A))}$, this follows from the fact that the algebra functor of $\lu{F}(A)$ is given on objects by $(X,\omega,(M_x)) \mapsto \int_{X} M_x d\omega$, and on morphisms of the form $(\la{Id}, (f(x))_{x \in X})$ (the morphisms between $(X,\omega,(M_x))$ and $(X,\omega,(N_x))$ involved in formation of the ultraproduct functor on $\lu{G}(\lu{F}(A))$) is given by $(\la{Id}, (f(x))_{x \in X}) \mapsto \int_X f(x) d\omega \circ \Delta_{\omega,\la{Id}}=\int_X f(x) d\omega$.

On the other hand, for a $T$-normal colax algebra $(B,m^B,a^B,i^B)$ we have an equality between $m^B$ and $m^{\lu{F}(\lu{G}(B))}$, this is simply because the ultrastructure of $\lu{G}(B)$ is given on objects by $\int_I M_i d \omega=m^B(I,\omega,(M_i)_{i \in I})$ and on morphisms by $ \int_I f_i d \omega=m^B(\la{Id}, (f_i)_{i \in I})$.

So the only thing left to check is the following: Starting from an ultracategory $A$, we should show that the Fubini transforms of $A$ and $\lu{G}\circ \lu{F}(A)$ agree, i.e. $\Delta^A= \Delta^{\lu{G}\circ \lu{F}(A)}$, same thing for the family of isomorphisms associated to principal ultrafilters $\epsilon$. On the other hand, starting from a $T$-normal colax algebra $B$, we should show that the colax associators of $B$ and  $\lu{F}\circ \lu{G}(B)$ agree, i.e. $a^{B} =a^{\lu{F}(\lu{G}(B))}$, same thing for the unitor, it needs to be checked that $i^{B}=i^{\lu{F}(\lu{G}(B))}$.

\begin{note*} \normalfont

In this section, we will show an isomorphism between the $1$-category of normal colax algebras with lax (left) algebra morphisms  and the $1$-category of ultracategories with left ultrafunctors. The next section will extend this equivalence to the respective $2$-categories.
\end{note*}
\subsection{Equality of the Fubini transforms}
\label{3}

Suppose we have an ultracategory $A$, we wish to show that $\Delta^A= \Delta^{\lu{G}\circ \lu{F}(A)}$. The main idea we use in the proof is the fact that $\int_S\nu_s d \omega =\int_{S \times T} \delta_{\pi(s,t)}d \int_S \iota_s\nu_s d \omega$. Here $\iota_s$ is the map from $T$ to $S \times T$ that sends $t$ to $(s,t)$, notice that for every $s \in S$, $\iota_s$ is a section to the projection map from $S \times T$ (we denote this projection map by $\pi$). Now let us have a  look at the following diagram:

\[\begin{tikzcd}
	{{\int_T M_t d(\int_S \nu_s  d \omega) =\int_T M_t d( \int_{S \times T} \delta_{\pi_{(s,t)}} d  \int_S \iota_s \nu_s d \omega)}} &&&& {\int_{S} \int_T M_t d \int_{S \times T} \delta_{\pi(s,t)}d\iota_s \nu_s d\omega} \\
	\\
	\\
	{\int_{S \times T}\int_T M_t d \delta_{\pi(s,t)}d \int_S \iota_s \nu_s d \omega} &&&& {\int_S \int_{S \times T} \int_TM_t d\delta_{\pi(s,t)} d \nu_s d \omega} \\
	\\
	\\
	{\int_{S \times T}M^{'}_{(s,t)} d \int_S \iota_s \nu_s d \omega } &&&& {\int_S \int_{S \times T} M^{'}_{(s,t)} d \iota_s \nu_s d \omega} &&&&& {} \\
	\\
	\\
	&&&& {\int_S \int_T M_t d \nu_s d \omega}
	\arrow["{{\Delta_{\omega,\nu_{\bullet}}}}", from=1-1, to=1-5]
	\arrow["{\Delta_{\int_S \iota_s \nu_s d\omega,\delta_{\pi\bullet}}}", from=1-1, to=4-1]
	\arrow["{\int_S \Delta_{\int_S\iota_s \nu_sd \omega,\delta_{\pi\bullet}}d\omega}"', from=1-5, to=4-5]
	\arrow["{{\int_S \Delta_{\iota_s\nu_s,\pi}d\omega}}"{description, pos=0.4}, curve={height=-85pt}, from=1-5, to=7-5]
	\arrow["{{\Delta_{\omega,\iota \nu_{\bullet}}}}", from=4-1, to=4-5]
	\arrow["{{\int_{S \times T} \epsilon_{T,\pi(s,t)}d \int_S \iota_s \nu_s d \omega}}", from=4-1, to=7-1]
	\arrow["{{\int_S \int_{S \times T} \epsilon_{T,\pi(s,t)}d  \iota_s \nu_s d \omega}}"', from=4-5, to=7-5]
	\arrow["{{\Delta_{\omega,\iota \nu_{\bullet}}}}", from=7-1, to=7-5]
	\arrow["{{\int_S \Delta_{\nu_s,\iota_s}d\omega}}"', from=7-5, to=10-5]
\end{tikzcd}\]
The Fubini transform associated to $A$ is just $\Delta_{\omega,\nu_{\bullet}}$, while on the other hand, the Fubini transform of $\lu{G}\circ \lu{F}(A)$ is given by the left composition in the diagram $\Delta_{\omega,\nu_{\bullet}}^{\lu{G}\circ \lu{F}(A)}= \int_S \Delta_{\nu_s,\iota_s}d\omega \circ \Delta_{\omega,\iota \nu_{\bullet}} \circ \int_{S \times T} \epsilon_{T,\pi(s,t)}d \int_S \iota_s \nu_s d \omega \circ \Delta_{\int_S \iota_s \nu_s d\omega,\delta_{\pi\bullet}}$. Recall that $M^{'}_{(s,t)}$ is defined by $M^{'}_{(s,t)}= M_{t}$ for every $(s,t)$, and $\pi$ is the projection map from $S\times T$ to $T$.

First notice that the upper square commutes by the third of Lurie's ultracategory axioms (axiom $C$), and the lower diagram is just the naturality of the categorical Fubini transform, so the only thing remaining is to show that the composition $\int_S \Delta_{\nu_s, \iota_s}d\omega \  \circ \  \int_S \int_{S \times T} \epsilon_{T,\pi(s,t)}d  \iota_s \nu_s d \omega \  \circ \int_S \Delta_{\int_S\iota_s \nu_sd \omega, \delta_{\pi\bullet}}d \omega$ is the identity, but notice that by definition \newline $\int_S \int_{S \times T} \epsilon_{T,\pi(s,t)}d  \iota_s \nu_s d \omega  \ \circ \int_S \Delta_{\int_S\iota_s \nu_sd \omega, \delta_{\pi\bullet}}d\omega$ is  $\int_S \Delta_{\iota_s\nu_s,\pi}d \omega$, then by using proposition $1.3.5$ of \cite{lurie2018ultracategories}, the composition of $\int_S \Delta_{\nu_s,\iota_s}d\omega$ and $\int_S \Delta_{ \iota_s\nu_s,\pi}d\omega$ is equal to $\int_S \Delta_{\nu_s, \la{Id}_{T}}d\omega$, which by corollary $1.3.6$ of \cite{lurie2018ultracategories} is the identity.

\subsection{Equality of the family of morphisms $\epsilon$ associated to principal ultrafilters}

Suppose that $A$ is an ultracategory we need to show that $\epsilon^A= \epsilon^{\lu{G}(\lu{F}(A))}$. In order to do so, we can look at the following diagram:

\[\begin{tikzcd}
	{\int_*\int_SM_{s}d\delta_{s_0}d*} \\
	\\
	{\int_S M_sd\delta_{s_0}} &&& {\int_*M_{s_0} d*} &&& {M_{s_0}}
	\arrow["{\int_*\epsilon_{S,s_0}d*}", from=1-1, to=3-4]
	\arrow["{\epsilon_{*,*}^{-1}=\Delta_{*,\delta_{s_0}}}", from=3-1, to=1-1]
	\arrow["{m_{* \mapsto s_0}}", from=3-1, to=3-4]
	\arrow["{\epsilon_{S,s_0}}", curve={height=30pt}, from=3-1, to=3-7]
	\arrow["{i_{M_{s_0}}=\epsilon_{*,*}}", from=3-4, to=3-7]
\end{tikzcd}\]

The upper composition corresponds to $\epsilon_{S,s_0}^{\lu{G}(\lu{F}(A))}$, which using the diagram above  equals $\epsilon_{S,s_0}^A$, by naturality of the $\epsilon_{*,*}$.

\subsection{Equality of the colax algebra associators}
On the other hand, suppose that we have a $T$-normal colax algebra $B$ with algebra functor $m$ (that means that there exists a $T^{'}$-colax algebra $m^{'}$ together with a colax algebra $2$-morphism $a^{'}$ and unitor $i^{'}$) .  Then we want to show the equality of colax associators of the $T^{'}$ colax algebra structures of $B$ and  $\lu{F}\circ \lu{G}(B)$. To show the equality of the colax algebra associators, we use the following diagram:

\[\begin{tikzcd}
	{\int_{\coprod_{S}T_s}M d \int_S\iota_s\gamma_s d\omega} &&&&&& {\int_{S \times \coprod_{S}T_s}M^{'}d  \int_{S}\iota_{(s',s)}\gamma_sd\omega} \\
	\\
	{\int_S \int_{T_s} M d \gamma_s d \omega} &&&&&& {\int_S \int_{\coprod_{S}T_s}M^{'} d \iota_s\gamma_sd\omega} \\
	\\
	&&&&&& {\int_S \int_{T_s}Md \gamma_s d\omega}
	\arrow["{(m^{'} \circ \mu_{B})_{(\la{Id},(\pi_s),(\la{Id})_{(s,t)})}=m^{'}_{(\la{Id},\pi,(\la{Id})_{s})}}", from=1-1, to=1-7]
	\arrow["{a^{'}}"', from=1-1, to=3-1]
	\arrow["{a^{'}}", from=1-7, to=3-7]
	\arrow["{(m^{'} \circ T^{'}m^{'})_{(\la{Id},(\pi_s),(\la{Id})_{(s,t)})}}"', from=3-1, to=3-7]
	\arrow["{\la{Id}}"', from=3-1, to=5-7]
	\arrow["{(m^{'} \circ T^{'}m^{'})_{(\la{Id},(\iota_s),(\la{Id})_{(s,t)})}}", from=3-7, to=5-7]
\end{tikzcd}\]
In this diagram, the ultraproduct is the ultraproduct of $\lu{G}(B)$, and $a^{'}$ on the left is \newline ${a^{'}}^{B}_{\ca{P}}$, where $\ca{P} = (S,\omega,(T_s)_{s\in S}, (M_{(s,t)})_{(s,t) \in \coprod_{s \in S}T_s},(\gamma_{s})_{s \in S})$, and the $a^{'}$ on the right is  ${a^{'}}^B_{\ca{Q}}$, where $\ca{Q}=(S,\omega,(S \times T_s)_{s \in S}, (M_{(s',s,t)})_{(s',s,t) \in S \times \coprod_{s \in S} T_s},(\iota_{(s',s)}\gamma_s)_{s \in S} )$.

Our goal is to show the equality of the two colax associators, which is equivalent to saying that the outermost diagram commutes, since ${a^{'}}^{\lu{F} (\lu{G}(B))}=\int_S \Delta^{\lu{G}(B)}_{\gamma_s, \iota_s}d\omega \ \circ \  \Delta^{\lu{G}(B)}_{\omega,\iota\gamma_{\bullet}} = ({m^{'}}^{B} \circ T^{'}{m^{'}}^{B})_{(\la{Id},(\iota_s),(\la{Id})_{(s,t)})} \  \circ \  {a^{'}}^B_{\ca{Q}} \  \circ {m^{'}}^B_{\la{Id},\pi,(\la{Id})_{s}}   $. Of course, one should justify why $\int_S \Delta^{\lu{G}(B)}_{\gamma_s, \iota_s}d\omega = ({m^{'}}^{B} \circ T^{'}{m^{'}}^{B})_{(\la{Id},(\iota_s),(\la{Id})_{(s,t)})}$, but this easily follows from the definition of $T^{'}$ (or $T$) on morphisms, and the definition of the ultraproduct functor of $\lu{G}(B)$.

Now, the upper diagram is commutative by naturality of $a^{'}$, and the lower diagram is commutative by the fact that the two morphisms in ${T^{'}}^{2}B$ compose to the identity, thus their image by $m^{'}  \circ T^{'}m^{'}$ is the identity.

\subsection{Equality of the unitors}
Now suppose that we have a $T$-normal colax algebra $B$, we want to show that $\lu{F}\circ \lu{G}(B)$ and $B$ have the same unitor. This follows from the fact that $m^{'}_{*\mapsto *}=\la{Id}$ (since functors preserve identity maps).

\subsection{Functoriality} \label{Functoriality}

Until now,  what we have shown is a bijective correspondence between ultrastructures on a category and normal colax algebra structures on the same category. Our next claim is that left ultrafunctors are lax (left) algebra morphisms, right ultrafunctors are colax (right) algebra morphisms, and ultrafunctors are algebra pseudo-morphisms. To do so, we present the following proof for the lax (left) algebra morphism case: \\  \noindent We are claiming that the data of a left ultrafunctor is the same as the data of a lax algebra morphism:

Recall that both are given by a natural family of morphisms $\sigma_{\omega}$ from $F(\int_X M_x d\omega)$ to $\int_XF(M_x)d\omega$ required to satisfy certain axioms (different axioms in  each case) (notice that the naturality condition in the case of monad algebras is a stronger condition than that in the case of ultracategory functors, so there is an additional thing to check in the case when we have a left ultrafunctor).

First, suppose that we have a lax algebra morphism. Let us denote by $m_1$ the algebra functor of the first category and by $m_2$ the algebra functor of the second one, then we want to show that its data satisfies the axioms $(1)$ and $(2)$ of the definition of left ultrafunctors of  \cite[Definition 1.4.1]{lurie2018ultracategories}. For axiom $(1)$, we use the following diagram:

\[\begin{tikzcd}
	{F(\int_I M_{i}d\delta_{i_0})} &&&& {\int_IF(M_{i})d\delta_{i_0}} \\
	&& {F(M_{i_0})} \\
	\\
	{F(\int_* M_{i_0}d*)} &&&& {\int_*F(M_{i_0})d*} \\
	&& {F(M_{i_0})}
	\arrow["{\sigma_{\delta_{i_0}}}"{description}, from=1-1, to=1-5]
	\arrow["{F(\epsilon_{I,i_{0}})}"{description}, from=1-1, to=2-3]
	\arrow[""{name=0, anchor=center, inner sep=0}, "{F(m_1)}"', from=1-1, to=4-1]
	\arrow[""{name=1, anchor=center, inner sep=0}, "{\epsilon_{I,i_0}}"{description}, from=1-5, to=2-3]
	\arrow[""{name=2, anchor=center, inner sep=0}, "{m_2}"', from=1-5, to=4-5]
	\arrow["1"', draw=none, from=2-3, to=4-1]
	\arrow["2"', draw=none, from=2-3, to=4-5]
	\arrow["{\la{Id}}"{description, pos=0.3}, shift left, from=2-3, to=5-3]
	\arrow["{\sigma_{*}}"{description}, from=4-1, to=4-5]
	\arrow["{F(\epsilon_{*,*})}"{description}, from=4-1, to=5-3]
	\arrow[""{name=3, anchor=center, inner sep=0}, "{\epsilon_{*,*}}"{description}, from=4-5, to=5-3]
	\arrow["5"', draw=none, from=1-1, to=1]
	\arrow["3"', draw=none, from=0, to=2]
	\arrow["4"', draw=none, from=4-1, to=3]
\end{tikzcd}\]

First, notice that squares $1,2$ (left and right front squares) commute, using naturality of $\epsilon_{*,*}$, the definition of $m$ (the algebra map) associated to an ultracategory and the fact that $\Delta_{*,\delta_{i_0}}$ from $\int_I M_i d \delta_{i_0}$ to $\int_* \int_{I}M_id\delta_{i_0} d*$ is the inverse of $\epsilon_{*,*}$, square $3$ (the back square of the triangular prism) and triangle $4$ commute by properties of the maps $\sigma$ (maps of left algebra morphism). Thus, since every arrow is invertible, we get that triangle $5$ commutes, which is exactly what we wanted to show.

Now, for axiom $(2)$ this can be deduced from the following diagram:

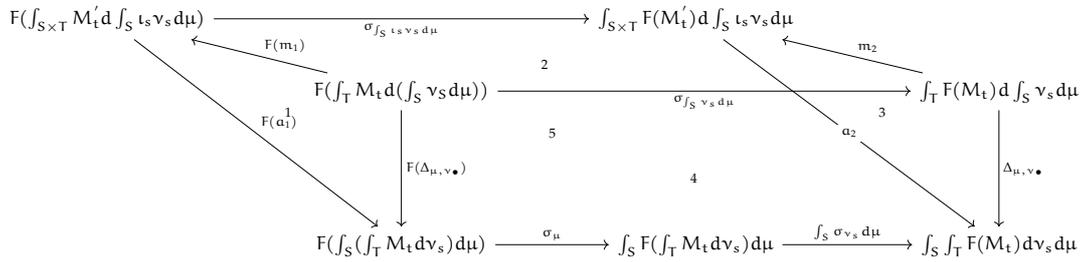
\[
\adjustbox{max width=\textwidth}{

\begin{tikzcd}
	& {F(\int_{S \times T}M^{'}_td\int_S\iota_s\nu_sd\omega)} &&&& {\int_{S\times T} F(M_t^{'} )d\int_S\iota_s\nu_sd\omega} & {} \\
	&&& {F(\int_T M_t d (\int_S \nu_s d \omega))} &&&& {\int_T F(M_t) d \int_S \nu_s d\omega} \\
	\\
	{} \\
	&&& {F(\int_S (\int_T M_t d\nu_s) d\omega)} && {\int_SF(\int_TM_td\nu_s)d \omega } && {\int_S \int_T F(M_t) d \nu_s d \omega} \\
	\\
	&&&&&& {}
	\arrow[""{name=0, anchor=center, inner sep=0}, "{\sigma_{\int_S \iota_s\nu_s d \omega}}"', from=1-2, to=1-6]
	\arrow[""{name=1, anchor=center, inner sep=0}, "{F(a_1)}"{description}, from=1-2, to=5-4]
	\arrow[""{name=2, anchor=center, inner sep=0}, "{a_2}"{description}, from=1-6, to=5-8]
	\arrow["{F(m_1)}"', from=2-4, to=1-2]
	\arrow[""{name=3, anchor=center, inner sep=0}, "{\sigma_{\int_S \nu_s d \omega}}"', from=2-4, to=2-8]
	\arrow["{F(\Delta_{\omega,\nu_{\bullet}})}", from=2-4, to=5-4]
	\arrow["4"', draw=none, from=2-4, to=5-8]
	\arrow["{m_2}"', from=2-8, to=1-6]
	\arrow["{\Delta_{\omega,\nu_{\bullet}}}", from=2-8, to=5-8]
	\arrow["{\sigma_\omega}", from=5-4, to=5-6]
	\arrow["{\int_S \sigma_{\nu_s}d\omega}", from=5-6, to=5-8]
	\arrow["2"', draw=none, from=0, to=3]
	\arrow["5"{description}, draw=none, from=1, to=2]
	\arrow["1"', draw=none, from=2-4, to=1]
	\arrow["3"', draw=none, from=2-8, to=2]
\end{tikzcd}
}
\]
We wish to show that square $4$ (front diagram) is commutative: \\ Square $5$ commutes by the fact that $F$ is a lax (left) morphism of colax algebras.\\
Triangles $1$ and $3$ commute by definition.\\
Square $2$ commutes by naturality.

\begin{note*} \normalfont
    The reader may object that we only required the commutativity of square $5$ in the case of ordinals in our definition. But we can extend the definition of the maps $\sigma_{\omega}$ in a unique way, which makes square $5$ commutative for every set.
\end{note*}

Now, on the other hand, suppose that we have a left ultrafunctor, then to show that this is a lax algebra morphism between the corresponding colax algebras: First, we check naturality in the family $(M_i)_{i \in I}$ (here the reader should note that the naturality condition imposed in the definition of  left ultrafunctors on the family $(\sigma_{\omega})$ is weaker than that imposed by the colax algebra definition, this is due to the fact that the ultracategory definition does not, in principle, relate the cases when we have different indexing sets (although it does by virtue of our work)):

In other words, we want to check the following:
\begin{theorem}
    \label{Stronger naturality}
Let $A$ and $B$ be ultracategories, and let $F$ be a left ultrafunctor between them. Suppose that we have $(X,f\omega,(M_x)_{x \in X})$ and $(Y,\iota_{Y^{'}}\omega,(N_y)_{y \in Y})$ and a morphism $( f:Y^{'}\subseteq Y \rightarrow X ,(\chi_y)_{y \in Y^{'}})$ of ${T}A$ between these two objects. Then the outermost diagram in \hyperref[diag6]{Diagram 6} of \hyperref[appendix2]{Appendix B} commutes:

\end{theorem}

\begin{proof}

Squares $1,7$ commute by axiom $2$ of \cite[definition 1.4.1]{lurie2018ultracategories}. \\
Squares $2,4,5$ commute by axiom $0$ (naturality). \\
Triangles $3,6$ commute by axiom $1$. \\
\end{proof}
Now, and in the same setting as in the previous argument, we need to check the following equality (there is another equality specific for the unit, but again, we omit it since the proof should be easy):

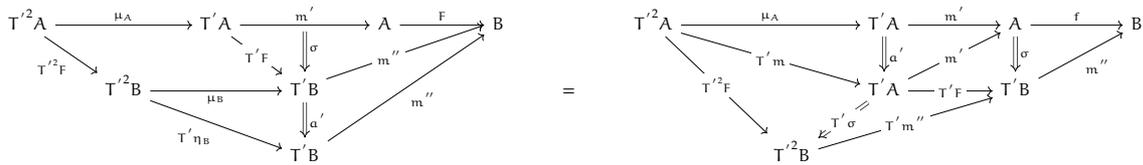
\[
\adjustbox{max width = \textwidth}{
\begin{tikzcd}
	{{T^{'}}^2A} && {T^{'}A} && A && B && {{T^{'}}^2A} &&& {T^{'}A} && A && B \\
	& {{T^{'}}^2B} && {T^{'}B} &&&& {=} &&&& {T^{'}A} && {T^{'}B} \\
	&&& {T^{'}B} &&&&&&& {{T^{'}}^2B}
	\arrow["{{\mu_A}}", from=1-1, to=1-3]
	\arrow["{{{T^{'}}^2F}}"', from=1-1, to=2-2]
	\arrow[""{name=0, anchor=center, inner sep=0}, "{{m^{'}}}", from=1-3, to=1-5]
	\arrow["{{T^{'}F}}"{description}, from=1-3, to=2-4]
	\arrow["F", from=1-5, to=1-7]
	\arrow["{{\mu_{A}}}", from=1-9, to=1-12]
	\arrow["{{T^{'}m^{'}}}"{description}, from=1-9, to=2-12]
	\arrow["{{{T^{'}}^2F}}"{description}, from=1-9, to=3-11]
	\arrow["{{m^{'}}}", from=1-12, to=1-14]
	\arrow["{{a^{'}}}", Rightarrow, from=1-12, to=2-12]
	\arrow["F", from=1-14, to=1-16]
	\arrow["\sigma", Rightarrow, from=1-14, to=2-14]
	\arrow["{{\mu_B}}"', from=2-2, to=2-4]
	\arrow["{{T^{'}m''}}"', from=2-2, to=3-4]
	\arrow["{{m^{''}}}"{description, pos=0.4}, shift left, from=2-4, to=1-7]
	\arrow["{{a^{''}}}", Rightarrow, from=2-4, to=3-4]
	\arrow["{{m^{'}}}"{description}, from=2-12, to=1-14]
	\arrow["{{T^{'}F}}"{description}, from=2-12, to=2-14]
	\arrow["{{T^{'}\sigma}}"{description}, Rightarrow, from=2-12, to=3-11]
	\arrow["{{m^{''}}}"', from=2-14, to=1-16]
	\arrow["{{m^{''}}}"', from=3-4, to=1-7]
	\arrow["{{T^{'}m^{''}}}"{description}, from=3-11, to=2-14]
	\arrow["\sigma", between={0.2}{1}, Rightarrow, from=0, to=2-4]
\end{tikzcd}
}
\]
\noindent Here $m^{'}$ and $m^{''}$ are the algebra functors for the strict monad $T^{'}$ structures, on the first and second category, respectively.

This translates to showing that the outermost diagram commutes:
\[
\adjustbox{max width =\textwidth}{

\begin{tikzcd}
	{F(\int_{\coprod_{s \in S}T_s}M_{(s,t)}d\int_S\iota_s\nu_sd\omega)} &&&& {\int_{\coprod_{s \in S}T_s} F(M_{(s,t)}) d\int_S\iota_s\nu_sd\omega} & {} \\
	\\
	&& {F(\int_S\int_{\coprod_{s \in S}T_s} M_{(s,t)} d \iota_s \nu_s d\omega)} && {\int_SF(\int_{\coprod_{s \in S}T_s} M_{(s,t)}) d \iota_s \nu_s d\omega} && {\int_S\int_{\coprod_{s \in S }T_s}F(M_{(s,t)})d\iota_s\nu_sd\omega} \\
	\\
	&& {F(\int_S (\int_{T_s} M_{(s,t)} d\nu_s) d\omega)} && {\int_SF(\int_{T_s}M_{(s,t)}d\nu_s)d \omega } && {\int_S \int_{T_s} F(M_{(s,t)}) d \nu_s d \omega} \\
	\\
	&&&&& {}
	\arrow["{{\sigma_{\int_S \iota_s\nu_s d \omega}}}"', from=1-1, to=1-5]
	\arrow[""{name=0, anchor=center, inner sep=0}, "{{F(\Delta_{\omega,\iota\nu_{\bullet}})}}", from=1-1, to=3-3]
	\arrow["2", draw=none, from=1-1, to=3-5]
	\arrow[""{name=1, anchor=center, inner sep=0}, "{{F(a^{'})}}"{description, pos=0.7}, from=1-1, to=5-3]
	\arrow[""{name=2, anchor=center, inner sep=0}, "{{\Delta_{\omega,\iota\nu_{\bullet}}}}", from=1-5, to=3-7]
	\arrow[""{name=3, anchor=center, inner sep=0}, "{a^{''}}"{description, pos=0.7}, from=1-5, to=5-7]
	\arrow["{{\sigma_{\omega}}}", from=3-3, to=3-5]
	\arrow["{{F(\int_S\Delta_{\nu_s,\iota_s}d\omega)}}"{description}, from=3-3, to=5-3]
	\arrow["{{\int_S \sigma_{\iota_s\nu_s}}}", from=3-5, to=3-7]
	\arrow["{{\int_SF(\Delta_{\nu_s,\iota_s})}}", from=3-5, to=5-5]
	\arrow["{{\int_S\Delta_{\nu_s,\iota_s}}}", from=3-7, to=5-7]
	\arrow["3", draw=none, from=5-3, to=3-5]
	\arrow["{{\sigma_\omega}}", from=5-3, to=5-5]
	\arrow[""{name=4, anchor=center, inner sep=0}, "{{\int_S \sigma_{\nu_s}d\omega}}", from=5-5, to=5-7]
	\arrow["1", draw=none, from=1, to=0]
	\arrow["5", draw=none, from=2, to=5-7]
	\arrow[draw=none, from=3-5, to=4]
	\arrow["4", draw=none, from=5-5, to=3]
\end{tikzcd}
}
\]
\noindent Triangles $1$ and $5$ commute by definition.\\
Square $2$ commutes by axiom $(2)$ of \cite[definition 1.4.1]{lurie2018ultracategories}. \\
Square $3$ commutes by naturality of $\sigma$ (axiom $0$).\\
Square $4$ commutes by the stronger version of naturality that we showed.\\
In this argument, $a^{'}$ denotes $a^{'}_{(S,\omega,(T_s)_{s \in S}, (\nu_s)_{s \in S},(M_{(s,t)})_{(s,t) \in \coprod_{s \in S} T_s})}$ in triangle $1$,\\and  $a^{''}$ denotes $a^{''}_{(S,\omega,(T_s)_{s \in S}, (\nu_s)_{s \in S},(F(M_{(s,t)}))_{(s,t) \in \coprod_{s \in S} T_s})}$ in triangle $5$, where $a^{'}$ and $a^{''}$ are the colax associators of $m^{'}$ and $m^{''}$ respectively.

\section{2-Categorical aspect} \label{natural transformations}

Up to this point, we can say that all our work was $1$-categorical, in other words, we have shown an equivalence (to be more precise, we have technically shown something stronger which is an isomorphism) between ultracategories with left (respectively right) ultrafunctors and normal colax algebras for this  monad with lax (respectively colax) algebra morphisms, this was done by showing that the structure of an ultracategory corresponds uniquely to a structure of a $T$-normal colax algebra. So at the $1$-category level we have shown a bijection between the additional data on a category needed to make a $T$-normal colax algebra/ultracategory, same on functors. But we have seen that ultracategories form a $2$-category together with either left or right ultrafunctors, and with the corresponding natural transformations of those, so to extend this isomorphism of categories to the level of the corresponding $2$-categories, we need to extend this bijection to the level of the data of $2$-morphisms.

We now claim that natural transformations of left/right ultrafunctors should correspond uniquely to the $2$-morphisms between left (lax)/right (colax) morphisms, more precisely, in both the definition of $2$-morphisms between lax (left) algebra morphisms and of natural transformations of left ultrafunctors, the data is that of a natural transformation that satisfies the commutativity of the following diagram:

\[\begin{tikzcd}
	{F(\int_I M_i d\omega)} &&& {G(\int_I M_i d\omega)} \\
	\\
	{\int_I F(M_i)d\omega} &&& {\int_I G(M_i)d\omega}
	\arrow["{\phi_{\omega}}", from=1-1, to=3-1]
	\arrow["{\psi_{\int_I M_i d\omega}}"', from=1-1, to=1-4]
	\arrow["{\phi^{'}_{\omega}}"', from=1-4, to=3-4]
	\arrow["{\int_I \psi_{M_i}d \omega}", from=3-1, to=3-4]
\end{tikzcd}\]

Or in the case of natural transformations of right ultrafunctors and $2$-morphisms between colax (right) algebra morphisms, the data is that of a natural transformation that satisfies the commutativity of the following diagram:

\[
\begin{tikzcd}
	{F(\int_I M_i d\omega)} &&& {G(\int_I M_i d\omega)} \\
	\\
	{\int_I F(M_i)d\omega} &&& {\int_I G(M_i)d\omega}
	\arrow["{\psi_{\int_I M_i d\omega}}"', from=1-1, to=1-4]
	\arrow["{\phi_{\omega}}"', from=3-1, to=1-1]
	\arrow["{\int_I \psi_{M_i}d \omega}", from=3-1, to=3-4]
	\arrow["{\phi^{'}_{\omega}}"', from=3-4, to=1-4]
\end{tikzcd}
\]
\section{What are pseudo-algebras for this monad} \label{algebras}
A final interesting question would be to ask for some description of the pseudo-algebras for this (pseudo)-monad. These are ultracategories for which the colax associator is invertible. We have stated that such objects have acquired the name of strong ultracategories. We state the following theorem:

\begin{theorem}
    Suppose that $A$ is the category of models of a first-order theory, then $A$ is a pseudo-algebra for the pseudo-monad $T$.
\end{theorem}
We are going to present the proof in the case where $A = \mathsf{Set}$, since the proof reduces to this case by Łoś theorem: Namely, suppose that we have a first-order signature $\lu{L}=\langle \go{S},\go{R},\go{F}\rangle$, and a family of axioms $\ca{A}$, then given a family of models of $\ca{A}$ $(M_i)_{i \in I}$, the ultraproduct of these models is constructed as ultraproducts of sets sortwise, that means that for every $S \in \go{S}$, we have $(\int_I M_i d\omega)^S := \int_I M_i^{S} d\omega$, and by Łoś theorem this ultraproduct (constructed as ultraproduct of sets sortwise) is a model of $\ca{A}$.

We want to show that the map $a: \int_{\coprod_S T_s}M_{(s,t)}d\int_S \iota_s \lambda_s d\omega \to \int_{S}\int_{T_s} M_{(s,t)} d \lambda_s d\omega $ by $(q_{(s,t)})\mapsto ((q_{(s,t)})) $ is a bijection (we assume that all sets are non-empty; if instead some of the $M_{(s,t)}$ are empty, one restricts attention to the subset of indices for which the relevant sets are non-empty: if this subset lies in the ultrafilter in question, the argument below applies  after restricting to it, and if it does not, then the corresponding ultraproduct is itself empty and the statement holds vacuously)(of course, the reader should verify that what was written is well-defined). Suppose that $((q_{s,t}))=((r_{s,t}))$, that means that for some subset $U \in \omega$, we have that for each $s \in U$, $(q_{(s,t)})= (r_{(s,t)})$, that means that for each $s$, there exists a set $V_s$ such that for each $t \in V_s$ we have $q_{(s,t)}=r_{(s,t)}$, now we simply take the set $\coprod_{s\in U} V_s \in \int_{S} \iota_s \lambda_s d \omega $ for which for every $(s,t)$ $q_{(s,t)}=r_{(s,t)}$, and we deduce that $(q_{(s,t)})=(r_{(s,t)})$, which shows injectivity, surjectivity is clear.

\begin{theorem}
    Suppose that $A$ is a category of models of continuous model theory, then $A$ is a pseudo-algebra for the pseudo-monad $T$.
\end{theorem}

Continuous model theory \cite{hart2023introduction} is a generalisation of classical model theory which allows an axiomatic treatment of a large class of metric structures, like Banach spaces, Hilbert spaces, etc. The ultrastructure of these categories was explored in \cite{Bundles}. Models of continuous model theory also include models of classical model theory. It can be shown that the ultracategory of models of continuous model theory is (ultra)-equivalent to the ultracategory of models of classical model theory.

Similarly to the case of discrete model theory, the proof reduces to checking this for the category $k\text{-}\mathsf{CompMet}$, the category of complete metric spaces for which the distance is bounded by $k$, recalling that in this category, the ultraproduct of a family of non-empty spaces $\int_I A_i d \omega$ is given by $\prod A_i / \sim$, where $\sim$ identifies two tuples if they are arbitrarily close on some set of the ultrafilter. Again, we want to show that the map $a:\int_{\coprod_S T_s}M_{(s,t)}d\int_S \iota_s \lambda_s d\omega \to \int_{S}\int_{T_s} M_{(s,t)} d \lambda_s d\omega $ by $(q_{(s,t)})\mapsto ((q_{(s,t)})) $ is a bijection, and again, we are going to assume that no metric space is empty to simplify the argument. Suppose that $((q_{(s,t)}))=((r_{(s,t)}))$. Let $\epsilon> 0$, there exists a set $S^{'} \in \omega$ such that for any $s \in S'$ $d((q_{(s,t)}),(r_{(s,t)})) < \epsilon $.
Now that means that for any $s \in S^{'}$, there exists a set $A_s \in \lambda_s$ such that for any $t \in A_s$ $d(q_{(s,t)}, r_{(s,t)} ) < \epsilon$. Now take the set $\coprod_{s \in S^{'}} A_s \in \int_{S} \iota_s \lambda_s d \omega$ then, for any $(s,t)$ in this set, we have $d(q_{(s,t)},r_{(s,t)})< \epsilon$, which implies that when passing to the ultraproduct $(q_{(s,t)})= (r_{(s,t)})$, this shows injectivity, again surjectivity is clear, and hence $a$ is an isomorphism.

\section*{Acknowledgment}
\label{sec:ack}
\addcontentsline{toc}{section}{\nameref{sec:ack}}
  I would like to acknowledge and thank my doctoral thesis supervisor Simon Henry for his valuable guidance, inputs, and feedback. This work was partially supported by \href{https://www.nserc-crsng.gc.ca/}{The Natural Sciences and Engineering Research Council of Canada (NSERC)}, funding reference number RGPIN-2020-067 awarded to Simon Henry, and by the Ontario Ministry of Colleges and Universities through \href{https://www.uottawa.ca/study/graduate-studies/funding-financing/awards/ontario}{the Ontario Graduate Scholarship and the QEII Graduate Scholarship in Science and Technology}.

  \newpage

\bibliographystyle{plainnat}
\bibliography{bib}

\newpage

\section*{Appendix A: Proof of Axiom $C$ for arbitrary sets}
\label{proof} \label{A}
\addcontentsline{toc}{section}{\nameref{proof}}

 We are going to work in the same setting as in \ref{axiom 3}  where we have arbitrary sets and not just ordinals, the reader may notice that the problem in generalising the argument from the case of ordinals to arbitrary sets is that the set $R \times S \times T$ could mean two different things either $(R \times S) \times T$  or $R \times (S \times T)$, to solve this let us have a look at \hyperref[diag7]{Diagram 7} of \hyperref[appendix2]{Appendix B}.

In \hyperref[diag7]{Diagram 7}, the argument for the commutativity of squares $1,2,3$ remains the same, but now we should justify the commutativity of square $4$, to do that let's have a look at \hyperref[diag8]{Diagram 8} of \hyperref[appendix2]{Appendix B}.
\\

We want to show that the front square commutes in \hyperref[diag8]{Diagram 8}, but we have already shown that the back diagram commutes in \ref{axiom 3}, and all side maps are isomorphisms, also squares $1,2,3,5$ are commutative by definition, it only remains to show that square $4$ is commutative towards this, let us have a look at the following diagram:

\adjustbox{max width =\textwidth}
{

\begin{tikzcd}[ampersand replacement=\&]
	\&\&\&\&\&\&\& {\int_{\phi(R)} \int_{\phi(S) \times \phi(T)}M^{'} d \int_S \iota_s \nu_{s} d \omega_r d\lambda} \\
	\\
	{\int_R \int_{S \times T}M^{'} d \int_S \iota_s \nu_{s} d \omega_r d\lambda} \&\&\&\&\&\& {\int_R \int_{\phi(S) \times \phi(T)}M^{'} d \int_S \iota_s \nu_{s} d \omega_r d\lambda} \\
	\\
	\&\&\& 1 \&\&\&\& {\int_{\phi(R)} \int_{\phi(S)} \int_{\phi(T)} M_t d\nu_s d \omega_r d\lambda} \\
	\\
	{\int_R \int_S \int_T M_t d\nu_s d \omega_r d\lambda} \&\&\&\&\&\& {\int_R \int_{\phi(S)} \int_{\phi(T)} M_t d\nu_s d \omega_r d\lambda}
	\arrow[""{name=0, anchor=center, inner sep=0}, "{{\int_{\phi(R)} a^{'}_{S ,\omega_r,(T^{'}),(\nu_s),(M^{'}_{(s,t)})}d\lambda}}"{description}, from=1-8, to=5-8]
	\arrow[""{name=1, anchor=center, inner sep=0}, "{{(m \circ Tm)_{(r\mapsto e_R(r)),((s,t) \mapsto (e_S(s),e_T(t)))}}}"{description}, from=3-1, to=1-8]
	\arrow["{{Tm_{(s,t) \mapsto (e_S(s), e_T(t))}= \int_R m_{(s,t) \mapsto (e_S(s), e_T(t))}}}"{description}, from=3-1, to=3-7]
	\arrow["{{\int_R a_{S ,\omega_r,(T^{'}),(\nu_s),(M^{'}_{(s,t)})}d\lambda}}"{description}, from=3-1, to=7-1]
	\arrow["{{m_{r\mapsto e_R(r)}}}"{description}, from=3-7, to=1-8]
	\arrow[""{name=2, anchor=center, inner sep=0}, "{{\int_R a^{'}_{S ,\omega_r,(T^{'}),(\nu_s),(M^{'}_{(s,t)})}d\lambda}}"{description}, from=3-7, to=7-7]
	\arrow["{{\int_R (m \circ Tm)_{s \mapsto e_S(s),t \mapsto e_T(t) }}}"{description}, from=7-1, to=7-7]
	\arrow["{{m_{r\mapsto e_R(r)}}}"{description}, from=7-7, to=5-8]
	\arrow["3"{description}, draw=none, from=1, to=3-7]
	\arrow["2"{description}, draw=none, from=2, to=0]
\end{tikzcd}
}

Squares 1-2 commute by definition, while triangle $3$ commutes by naturality of $m$.

\section*{Appendix B: Large diagrams}
\label{appendix2}
\addcontentsline{toc}{section}{\nameref{appendix2}}
\pagenumbering{gobble}
\begin{figure}

\caption*{Diagram 1:}
\[
\rotatebox{90}{
\adjustbox{max width = 70 em, max height = 70 em }
{
\begin{tikzcd}[ampersand replacement=\&]
	{\int_{\coprod_{s \in S}T_s}M_{(s,t)}d \int_S \iota_s \nu_s d f\omega } \&\& {\int_S\int_{\coprod_{s \in S}T_s} M_{(s,t)} d \iota_s \nu_s df\omega} \&\& {\int_S\int_{T_s} M_{(s,t)} d \nu_s df\omega} \\
	\&\&\& 2 \\
	\&\& {\int_{Y^{'}} \int_{\coprod_{s \in S}T_s}M_{(s,t)}d \iota_{f(y)}g_{y}\gamma_{y}d\omega} \&\& {\int_{Y^{'}}\int_{T_{f(y)}}M_{(f(y),t)}dg_y\gamma_y d\omega} \\
	{\int_{\coprod_{y \in Y^{'}}X^{'}_y}M_{(f(y), g_y(x))}d\int_{Y^{'}}\iota_{y} \gamma_{y}d\omega} \&\& {\int_{Y^{'}}\int_{\coprod_{y \in Y^{'}}X^{'}_y}M_{(f(y), g_y(x))}d\iota_{y} \gamma_{y}d\omega} \&\& {\int_{Y^{'}} \int_{X^{'}_y} M_{(f(y),g_y(x))} d \gamma_yd\omega} \\
	\\
	\\
	{\int_{\coprod_{y \in Y^{'}}X^{'}_y}N_{(y,x)}d\int_{Y^{'}}\iota_{y} \gamma_{y}d\omega} \&\& {\int_{Y^{'}}\int_{\coprod_{y \in Y^{'}} X^{'}_y}N_{(y,x)}d\iota_{y} \gamma_{y}d\omega} \&\& {\int_{Y^{'}} \int_{X^{'}_y} N_{(y,x)}d \gamma_yd\omega} \\
	\&\& {\int_{Y}\int_{\coprod_{y \in Y^{'}} X^{'}_y}N_{(y,x)}d\iota_{y} \gamma_{y}d\omega^{'}} \&\& {\int_{Y} \int_{X^{'}_y} N_{(y,x)}d \gamma_yd\omega^{'}} \\
	\\
	{\int_{\coprod_{y \in Y}X_y}N_{(y,x)} d \int_{Y} \iota_y \gamma_{y}^{'}d\omega^{'}} \&\& {\int_{Y}\int_{\coprod_{y \in Y }X_y}N_{(y,x)}d\iota_{y} \gamma^{'}_{y}d\omega^{'}} \&\& {\int_{Y} \int_{X_y} N_{(y,x)}d \gamma^{'}_yd\omega^{'}}
	\arrow["{{\Delta_{f\omega,\iota \nu_{\bullet}}}}", from=1-1, to=1-3]
	\arrow["{\Delta_{\int_{Y^{'}} \iota_y \gamma_y d \omega , \bar{f}}}"', from=1-1, to=4-1]
	\arrow["1"', draw=none, from=1-1, to=4-3]
	\arrow["{{\int_S \Delta_{\nu_s,\iota_s}}df\omega}", from=1-3, to=1-5]
	\arrow["{\Delta_{\omega,f}}"{description}, from=1-3, to=3-3]
	\arrow["{\Delta_{\omega,f}}", from=1-5, to=3-5]
	\arrow["{\int_{Y^{'}} \Delta_{\gamma_y, \iota_{f(y)}} d\omega}", from=3-3, to=3-5]
	\arrow["{\int_{Y^{'}}\Delta_{ \iota_{y}\gamma_{y},\bar{f}} d\omega}"{description}, from=3-3, to=4-3]
	\arrow["3"{description}, draw=none, from=3-3, to=4-5]
	\arrow["{\int_{Y^{'}}\Delta_{\gamma_{y},g_y} d \omega}", from=3-5, to=4-5]
	\arrow["{{\Delta_{\omega,\iota \gamma_{\bullet}}}}"', from=4-1, to=4-3]
	\arrow["{\int_{\coprod_{y \in Y^{'}}X^{'}_y}\psi_{(y,x)}d\int_{Y^{'}}\iota_{y} \gamma_{y}d\omega}"', from=4-1, to=7-1]
	\arrow["4", draw=none, from=4-1, to=7-3]
	\arrow["{\int_{Y^{'}}\Delta_{\gamma_y,\iota_y} d\omega}"', from=4-3, to=4-5]
	\arrow["{\int_{Y^{'}}\int_{\coprod_{y \in Y^{'}}X^{'}_y}\psi_{(y,x)}d \iota_y\gamma_{f(y)}d\omega}"{description}, from=4-3, to=7-3]
	\arrow["5", draw=none, from=4-3, to=7-5]
	\arrow["{\int_{Y^{'}}\int_{X_y}\psi_{(y,x)}d\iota_y\gamma_{f(y)}d\omega}", from=4-5, to=7-5]
	\arrow["{{\Delta_{\omega,\iota \gamma_{\bullet}}}}", from=7-1, to=7-3]
	\arrow["{\Delta^{-1}_{\int_{Y}\iota_y\gamma_yd\omega, j}}"', from=7-1, to=10-1]
	\arrow["6", draw=none, from=7-1, to=10-3]
	\arrow["{\int_{Y^{'}}  \Delta_{\gamma_{y}, \iota_{y}} d\omega}", from=7-3, to=7-5]
	\arrow["{\Delta^{-1}_{\omega,i}}"{description}, from=7-3, to=8-3]
	\arrow["7", draw=none, from=7-3, to=8-5]
	\arrow["{\Delta^{-1}_{\omega,i}}", from=7-5, to=8-5]
	\arrow["{\int_{Y} \Delta_{\gamma_y,\iota_{y}}d \omega^{'}}"{description}, from=8-3, to=8-5]
	\arrow["{\int_{Y}  \Delta^{-1}_{\iota_y \gamma_y,j} d\omega^{'}}"{description}, from=8-3, to=10-3]
	\arrow["8", draw=none, from=8-3, to=10-5]
	\arrow["{\int_{Y} \Delta^{-1}_{\iota_y,\gamma_y^{'}} d \omega'}", from=8-5, to=10-5]
	\arrow["{{\Delta_{\omega^{'},\iota \gamma'_{\bullet}}}}", from=10-1, to=10-3]
	\arrow["{\int_{Y}  \Delta_{\gamma_{y}^{'}, \iota_{y}} d\omega'}", from=10-3, to=10-5]
\end{tikzcd}
\label{diag1}
}
}
\]
\end{figure}

\begin{figure}
\caption*{Diagram 2:}
\[
\rotatebox{90}{
\adjustbox{max width = 70 em, max height = 70 em }{ 
\begin{tikzcd}[ampersand replacement=\&]
	{\int_{\coprod_{R}\coprod_{S_r}T_{(r,s)}}M d \rho} \&\&\& {\int_{\coprod_{R}S_r} \int_{\coprod_{R}\coprod_{S_r}T_{(r,s)}} Md \iota_{(r,s)} \nu_{(r,s)} d \int_R \iota_r \omega_r d\lambda} \&\&\&\& {\int_{\coprod_R S_r} \int_{T_{(r,s)}}M  d\nu_{s}d \int_R \iota_r\omega_r d \lambda} \\
	\\
	\&\&\&\&\&\&\& {\int_{R} \int_{\coprod_{R}S_r} \int_{T_{(r,s)}} M d \nu_s d \iota_r \omega_r d \lambda} \\
	{\int_{R} \int_{\coprod_{R}\coprod_{S_r}T_{(r,s)}} Md \iota_r \int_{S_r} \iota_{s}\nu_{(r,s)}d \omega_{r}d\lambda} \\
	\\
	{\int_R \int_{\coprod_{S_r}T_{(r,s)}}M d \int_{S_r} \iota_s \nu_{(r,s)} d \omega_r d\lambda} \&\&\& {\int_R \int_{S_r} \int_{\coprod_{S_r}T_{(r,s)}} Md \iota_{s} \nu_{(r,s)} d \omega_r d \lambda} \&\&\&\& {\int_R \int_{S_r} \int_{T_{(r,s)}} M d\nu_{(r,s)} d \omega_r d\lambda}
	\arrow["{{\Delta_{\int_R\iota_r \omega_r d\lambda,\iota \nu_{\bullet}}}}", from=1-1, to=1-4]
	\arrow["{{\Delta_{\lambda,\iota_{\bullet}\int_{S_{\bullet}\iota_s \nu_{(\bullet,s)}d\omega_{\bullet}}}}}", from=1-1, to=4-1]
	\arrow["{{\int_{\coprod_{R}S_r} \Delta_{\omega_r,\iota_{(r,s)}}}}", from=1-4, to=1-8]
	\arrow["{{\Delta_{\lambda,\iota\omega_{\bullet}}}}", from=1-8, to=3-8]
	\arrow["{{\int_R \Delta_{\omega_r,\iota_r}d\lambda}}", from=3-8, to=6-8]
	\arrow["{{\int_R \Delta_{\int_S \iota_s \nu_{(r,s)} d \omega_r,\iota_r}d\lambda}}", from=4-1, to=6-1]
	\arrow["{{\int_R \Delta_{\omega_r,\int_S \iota_s \nu_{(r,\bullet)} }d\lambda}}", from=6-1, to=6-4]
	\arrow["{{\int_R \int_{S_r} \Delta_{\nu_{(r,s)},\iota_s}d\omega_rd\lambda}}", from=6-4, to=6-8]
\end{tikzcd}
\label{diag2}
}
}
\]
\end{figure}

\begin{figure}
\caption*{Diagram 3:}
\[
\rotatebox{90}{
\adjustbox{max width = 70em, max height = 70 em }{
\begin{tikzcd}[ampersand replacement=\&]
	{} \&\&\&\&\&\&\&\&\& \\
	\\
	\\
	\\
	\\
	\\
	{\int_{\coprod_{R}\coprod_{S_r}T_{(r,s)}}M d \rho} \&\&\& {\int_{\coprod_{R}S_r} \int_{\coprod_{R}\coprod_{S_r}T_{(r,s)}} Md \iota_{(r,s)} \nu_{(r,s)} d \int_R \iota_r \omega_r d\lambda} \&\&\&\&\&\& {\int_{\coprod_R S_r} \int_{T_{(r,s)}}M  d\nu_{s}d \int_R \iota_r\omega_r d \lambda} \\
	\\
	\&\&\& {\int_R(\int_{\coprod_{R}S_r}(\int_{\coprod_{R}\coprod_{S_r}T_{(r,s)}} Md \iota_r\iota_s\nu_{(r,s)}) d  \iota_r \omega_r) d \lambda} \&\&\&\&\&\& {\int_{R} \int_{\coprod_{R}S_{r}} \int_{T_{(r,s)}} M d \nu_s d \iota_r \omega_r d \lambda} \\
	{\int_{R} \int_{\coprod_{R}\coprod_{S_r}T_{(r,s)}} Md \iota_r \int_{S_r} \iota_{s}\nu_{(r,s)}d \omega_{r}d\lambda} \&\&\& {\int_R\int_{S_r} \int_{\coprod_{R}\coprod_{S_r}T_{(r,s)}} Md \iota_r \iota_s\nu_{(r,s)}d\omega_r d \lambda} \\
	\\
	{\int_R \int_{\coprod_{S_r}T_{(r,s)}}M d \int_{S_r} \iota_s \nu_{(r,s)} d \omega_r d\lambda} \&\&\& {\int_R \int_{S_r} \int_{\coprod_{S_r}T_{(r,s)}} Md \iota_{s} \nu_{(r,s)} d \omega_r d \lambda} \&\&\&\&\&\& {\int_R \int_{S_r} \int_{T_{(r,s)}} M d\nu_{(r,s)} d \omega_r d\lambda}
	\arrow["{{{\Delta_{\int_R\iota_r \omega_r d\lambda,\iota \nu_{\bullet}}}}}", from=7-1, to=7-4]
	\arrow["1", draw=none, from=7-1, to=9-4]
	\arrow["{{\Delta_{\lambda,\iota_{\bullet}\int_{S_{\bullet}\iota_s \nu_{(\bullet,s)}d\omega_{\bullet}}}}}", from=7-1, to=10-1]
	\arrow["{{{\int_{\coprod_{R}S_r} \Delta_{\omega_r,\iota_{(r,s)}}}}}", from=7-4, to=7-10]
	\arrow["{{{\Delta_{\lambda,\iota\omega_{\bullet}}}}}", from=7-4, to=9-4]
	\arrow["3", draw=none, from=7-4, to=9-10]
	\arrow["{{{\Delta_{\lambda,\iota\omega_{\bullet}}}}}", from=7-10, to=9-10]
	\arrow["{{{\int_{R}\int_{\coprod_{S_r}T_{(r,s)}}\Delta_{\nu_{(r,s)},\iota_r \iota_s} d \int_R\iota_r \omega_r d \lambda}}}", from=9-4, to=9-10]
	\arrow["{{{\int_R \Delta_{\omega_r,\iota_r}d \lambda}}}", from=9-4, to=10-4]
	\arrow[""{name=0, anchor=center, inner sep=0}, "{{{\int_R \Delta_{\omega_r,\iota_r}d\lambda}}}", from=9-10, to=12-10]
	\arrow["{{{\int_R\Delta_{\iota_r\omega_r,\iota\nu_{\bullet}}}}}", from=10-1, to=9-4]
	\arrow[""{name=1, anchor=center, inner sep=0}, "{{{\int_R\Delta_{\omega_r,\iota\nu_{\bullet}}}}}", from=10-1, to=10-4]
	\arrow["{{\int_R \Delta_{\int_S \iota_s \nu_{(r,s)} d \omega_r,\iota_r}d\lambda}}", from=10-1, to=12-1]
	\arrow["5", draw=none, from=10-1, to=12-4]
	\arrow["{{{\int_R \int_S \Delta_{\iota_s\nu_{(r,s)},\iota_r}}}}", from=10-4, to=12-4]
	\arrow["{{{\int_{R}\int_S\Delta_{\nu_{(r,s)},\iota_{(r,s)}}d\omega_rd\lambda}}}", from=10-4, to=12-10]
	\arrow["{{{\int_R \Delta_{\omega_r,\int_S \iota_s \nu_{\bullet} }d\lambda}}}", from=12-1, to=12-4]
	\arrow[""{name=2, anchor=center, inner sep=0}, "{{{\int_R \int_{S} \Delta_{\nu_{(r,s)},\iota_s}d\omega_rd\lambda}}}", from=12-4, to=12-10]
	\arrow["2"{description}, draw=none, from=9-4, to=1]
	\arrow["4"{description}, draw=none, from=9-4, to=0]
	\arrow["6"{marking, allow upside down}, draw=none, from=10-4, to=2]
\end{tikzcd}
}
}
\]
\label{diag3}
\end{figure}

\begin{figure}
\caption*{Diagram 4:}
\[
\rotatebox{90}{

\adjustbox{max width = 70 em, max height = 70 em }{

\begin{tikzcd}[ampersand replacement=\&]
	{\int_TM_td\rho} \&\&\& {\int_{R \times T} M^{'}_{(r,t)} d \int_R \iota_r \int_S \nu_s d \omega_r d \lambda} \&\&\& {\int_R \int_T M_t d \int_S \nu_s d \omega_r d\lambda} \\
	\\
	{\int_{S \times T} M^{'}_{(s,t)} d \int_S \iota_s \nu_s d(\int_R\omega_r d\lambda)} \&\&\&\&\&\& {\int_R \int_{S \times T}M^{'} d \int_S \iota_s \nu_s d \omega_r d\lambda} \\
	\\
	{\int_S (\int_T M_t d\nu_s)  d\int_R \omega_r d \lambda} \&\&\& {\int_{R \times S} \int_{T}M_t  d\nu_{s}d \int_R \iota_r\omega_r d \lambda} \&\&\& {\int_R \int_S \int_T M_t d\nu_s d \omega_r d\lambda}
	\arrow["{{m_{(r,t) \mapsto t}}}", from=1-1, to=1-4]
	\arrow["{{m_{(s,t) \mapsto t }}}"', from=1-1, to=3-1]
	\arrow["{{a_{R,\lambda,(T^{'}),(\iota_r\int_S\nu_s d\omega_r),(M^{'})}}}", from=1-4, to=1-7]
	\arrow["{{\int_R m_{(s,t) \mapsto t} d\lambda}}", from=1-7, to=3-7]
	\arrow["{{a_{S,\int_R\omega_rd\lambda,(T^{'}_s),(\nu_s),(M^{'}_{(s,t)}) }}}"', from=3-1, to=5-1]
	\arrow["{{\int_R a_{S ,\omega_r,(T^{'}),(\nu_s),(M^{'}_{(s,t)})}d\lambda}}", from=3-7, to=5-7]
	\arrow["{{m_{(r,s)\mapsto s}}}", from=5-1, to=5-4]
	\arrow["{{a_{R,\lambda,(S^{'}),(\omega_r),(\int_TM_td\nu_s)^{'}}}}", from=5-4, to=5-7]
\end{tikzcd}
}
}
\]
\label{diag4}
\end{figure}

\begin{figure}
\caption*{Diagram 5:}
\[
\rotatebox{90}{

\adjustbox{ scale={1}{1} }{
\begin{tikzcd}[ampersand replacement=\&]
	{\int_TM_td\rho} \&\&\& {\int_{R \times T} M^{'}_{(r,t)} d \int_R \iota_r \int_S \nu_s d \omega_r d \lambda} \&\&\&\& {\int_R \int_T M_t d \int_S \nu_s d \omega_r d\lambda} \\
	\\
	{\int_{S \times T} M^{'}_{(s,t)} d \int_S \iota_s \nu_s d(\int_R\omega_r d\lambda)} \&\&\& {\int_{R \times S  \times T}M^{''}_{(r,s,t)}d \int_R \iota_r (\int_S \iota_s \nu_s d\omega_r) d\lambda} \&\&\&\& {\int_R \int_{S \times T}M^{'} d \int_S \iota_s \nu_s d \omega_r d\lambda} \\
	\\
	{\int_S (\int_T M_t d\nu_s)  d\int_R \omega_r d \lambda} \&\&\& {\int_{R \times S} \int_{T}M_t  d\nu_{s}d \int_R \iota_r\omega_r d \lambda} \&\&\&\& {\int_R \int_S \int_T M_t d\nu_s d \omega_r d\lambda}
	\arrow["{{{{m^{'}_{(r,t) \mapsto t}}}}}", from=1-1, to=1-4]
	\arrow["{{{{m^{'}_{(s,t) \mapsto t }}}}}"', from=1-1, to=3-1]
	\arrow["1", draw=none, from=1-1, to=3-4]
	\arrow[""{name=0, anchor=center, inner sep=0}, "{{{{a_{R,\lambda,(T^{'}),(\iota_r\int_S\nu_s d\omega_r),(M^{'})}}}}}", from=1-4, to=1-8]
	\arrow["{{{{m^{'}_{(r,s,t)\mapsto (r,t)}}}}}", from=1-4, to=3-4]
	\arrow["{{{{\int_R m^{'}_{(s,t) \mapsto t} d\lambda}}}}", from=1-8, to=3-8]
	\arrow["{{{{m^{'}_{(r,s,t)\mapsto (s,t)}}}}}", from=3-1, to=3-4]
	\arrow["{{{{a^{'}_{S,\int_R\omega_rd\lambda,(T^{'}_s),(\nu_s),(M^{'}_{(s,t)}) }}}}}"', from=3-1, to=5-1]
	\arrow["3", draw=none, from=3-1, to=5-4]
	\arrow[""{name=1, anchor=center, inner sep=0}, "{{{{a^{'}_{R,\lambda,(S\times T)^{'},(\iota_r\int_S \iota_s\nu_sd\omega_r),(M^{'})}}}}}", from=3-4, to=3-8]
	\arrow["{{{{a^{'}_{R\times S, \int_R \iota_r\omega_r d \lambda,(T^{''}_{(r,s)}),(\iota_s\nu_{s}),(M^{''}_{(r,s,t)})}}}}}", from=3-4, to=5-4]
	\arrow["{{{{\int_R a^{'}_{S ,\omega_r,(T^{'}),(\nu_s),(M^{'}_{(s,t)})}d\lambda}}}}", from=3-8, to=5-8]
	\arrow["{{{{m^{'}_{(r,s)\mapsto s}}}}}", from=5-1, to=5-4]
	\arrow[""{name=2, anchor=center, inner sep=0}, "{{{{a^{'}_{R,\lambda,(S^{'}),(\omega_r),(\int_TM_td\nu_s)^{'}}}}}}", from=5-4, to=5-8]
	\arrow["2"{description}, draw=none, from=0, to=1]
	\arrow["4"{description, pos=0.4}, draw=none, from=1, to=2]
\end{tikzcd}
}
}
\]
\label{diag5}
\end{figure}

\begin{figure}
\caption*{Diagram 6:}
\[
\rotatebox{90}
{
\adjustbox{max width= 70em,max height = 70em}
{
\begin{tikzcd}[ampersand replacement=\&]
	{F(\int_X M_x d f \omega)} \&\& {F(\int_{Y^{'}}\int_X M_x d \delta_{f(y)}d\omega)} \&\& {F(\int_{Y^{'}} M_{f(y)} d  \omega)} \&\& { F(\int_{Y^{'}} N_y d \omega)} \&\& {F(\int_{Y^{'}}\int_{Y}N_yd\delta_{\iota_{Y^{'}}(y)}d\omega)} \&\& {F(\int_Y N_yd\iota_{Y^{'}}\omega)} \\
	\&\& {\int_{Y^{'}}F(\int_XM_xd \delta_{f(y)})d\omega} \&\&\&\&\&\& {\int_{Y^{'}}F(\int_YN_yd\delta_{\iota_{Y^{'}}(y)})d\omega} \\
	{\int_XF(M_x)dfd\omega} \&\& {\int_{Y^{'}}\int_{X}F(M_x)d \delta_{f(y)}d\omega} \&\& {\int_{Y^{'}} F(M_{f(y)}) d  \omega} \&\& { \int_{Y^{'}} F(N_y) d \omega} \&\& {\int_{Y^{'}}\int_{Y}F(N_y)d\delta_{\iota_{Y^{'}}(y)}d\omega} \&\& {\int_{Y}F(N_y)d\iota_{Y^{'}}\omega}
	\arrow["{{F(\Delta_{\omega,\delta_{f(\bullet)}})}}"', from=1-1, to=1-3]
	\arrow["{{\sigma_{\iota_{Y^{'}}\omega}}}", from=1-1, to=3-1]
	\arrow["1", draw=none, from=1-1, to=3-3]
	\arrow["{{F(\int_{Y^{'}}\epsilon d\omega)}}"', from=1-3, to=1-5]
	\arrow["{{\sigma_{\omega}}}"', from=1-3, to=2-3]
	\arrow["2", draw=none, from=1-3, to=3-5]
	\arrow["{{F(\int_{Y^{'}} \chi_y d\omega)}}", from=1-5, to=1-7]
	\arrow["{{\sigma_{\omega}}}"', from=1-5, to=3-5]
	\arrow["4"{description}, draw=none, from=1-5, to=3-7]
	\arrow["{{F(\int_{Y^{'}}\epsilon^{-1})}}", from=1-7, to=1-9]
	\arrow["{{\sigma_{\omega}}}", from=1-7, to=3-7]
	\arrow["{{F(\Delta_{\omega,\delta_{\iota_{Y^{'}\bullet}}})^{-1}}}", from=1-9, to=1-11]
	\arrow["{{\sigma_{\omega}}}", from=1-9, to=2-9]
	\arrow["7"{description}, draw=none, from=1-11, to=3-9]
	\arrow["{{\sigma_{\omega}}}", from=1-11, to=3-11]
	\arrow["{{\int_{Y^{'}}\sigma_{\delta_{f(y)}}}d\omega}"', from=2-3, to=3-3]
	\arrow["{{\int_{Y^{'}}F(\epsilon)}d\omega}"{pos=0.6}, from=2-3, to=3-5]
	\arrow["{{\int_{Y^{'}}\sigma_{\delta_{\iota_{Y^{'}}(y)}}}}", from=2-9, to=3-9]
	\arrow["{{\Delta_{\omega,\delta_{f(\bullet)}}}}", from=3-1, to=3-3]
	\arrow[""{name=0, anchor=center, inner sep=0}, "{{\int_{Y^{'}}\epsilon d\omega}}"', from=3-3, to=3-5]
	\arrow["{{\int_{Y^{'}}F(\chi_yd\omega)}}"', from=3-5, to=3-7]
	\arrow[""{name=1, anchor=center, inner sep=0}, "{{\int_{Y^{'}}\epsilon^{-1}}d\omega}", from=3-7, to=2-9]
	\arrow[""{name=2, anchor=center, inner sep=0}, "{{\int_{Y^{'}}F(\epsilon^{-1}})d\omega}"', from=3-7, to=3-9]
	\arrow["{\Delta_{\omega,\delta_{\iota_{Y^{'}\bullet}}}}"', from=3-9, to=3-11]
	\arrow["5", draw=none, from=1-7, to=1]
	\arrow["3"{description}, draw=none, from=2-3, to=0]
	\arrow["6"{description}, draw=none, from=2-9, to=2]
\end{tikzcd}
}
}
\]
\label{diag6}
\end{figure}

\begin{figure}
\caption*{Diagram 7:}
\centering
\rotatebox{90}
{
\adjustbox{scale={0.67}{0.67}}
{
\begin{minipage}{40cm}
\[
\begin{tikzcd}[ampersand replacement=\&,row sep=30, column sep= 20]
    {\int_TM_td\rho} \&\&\&\& {\int_{R \times T} M^{'}_{(r,t)} d \int_R \iota_r \int_S \nu_s d \omega_r d \lambda} \&\&\& {\int_R \int_T M_t d \int_S \nu_s d \omega_r d\lambda} \\
    \\
    {\int_{S \times T} M^{'}_{(s,t)} d \int_S \iota_s \nu_s d(\int_R\omega_r d\lambda)} \&\& {\int_{(R \times S) \times T}M^{''}_{(r,s,t)}d \int_R \iota_r (\int_S \iota_s \nu_s d\omega_r) d\lambda} \&\& {\int_{R \times (S \times T)}M^{''}_{(r,s,t)}d \int_R \iota_r (\int_S \iota_s \nu_s d\omega_r) d\lambda} \&\&\& {\int_R \int_{S \times T}M^{'} d \int_S \iota_s \nu_{s} d \omega_r d\lambda} \\
    \\
    {\int_S (\int_T M_t d\nu_s)  d\int_R \omega_r d \lambda} \&\& {\int_{R \times S} \int_{T}M_t  d\nu_{s}d \int_R \iota_r\omega_r d \lambda} \&\&\&\&\& {\int_R \int_S \int_T M_t d\nu_s d \omega_r d\lambda}
    \arrow["{\mathfrak{a}}", from=1-1, to=1-5]
    \arrow["{\mathfrak{b}}"', from=1-1, to=3-1]
    \arrow[""{name=0, anchor=center, inner sep=0}, "{\mathfrak{c}}", from=1-5, to=1-8]
    \arrow[""{name=2, anchor=center, inner sep=0}, "{\mathfrak{d}}"{description}, from=1-5, to=3-3]
    \arrow[""{name=3, anchor=center, inner sep=0}, "{\mathfrak{e}}", from=1-5, to=3-5]
    \arrow["{\mathfrak{f}}", from=1-8, to=3-8]
    \arrow["{\mathfrak{g}}", from=3-1, to=3-3]
    \arrow["{\mathfrak{h}}"', from=3-1, to=5-1]
    \arrow["{\mathfrak{i}}", from=3-3, to=3-5]
    \arrow["{\mathfrak{j}}"{description}, from=3-3, to=5-3]
    \arrow[""{name=4, anchor=center, inner sep=0}, "{\mathfrak{k}}", from=3-5, to=3-8]
    \arrow["{\mathfrak{l}}", from=3-8, to=5-8]
    \arrow[""{name=5, anchor=center, inner sep=0}, "{\mathfrak{m}}", from=5-1, to=5-3]
    \arrow[""{name=6, anchor=center, inner sep=0}, "{\mathfrak{n}}", from=5-3, to=5-8]
    \arrow["1", draw=none, from=0, to=3-3]
    \arrow["5", draw=none, from=3, to=2]
    \arrow["2", draw=none, from=1, to=4]
    \arrow["3", draw=none, from=3-1, to=5]
    \arrow["4", draw=none, from=3-5, to=6]
\end{tikzcd}
\]

\vspace{2em}
\LARGE
\begin{align*} 
 \mathfrak{a} &= m_{(r,t) \mapsto t} \\
 \mathfrak{b} &= m_{(s,t) \mapsto t } \\
 \mathfrak{c} &= a_{R,\lambda,(T^{'}_s),(\iota_r\int_S\nu_s d\omega_r),(M^{'})} \\
 \mathfrak{d} &= m_{((r,s),t)\mapsto (r,t)} \\
 \mathfrak{e} &= m_{(r,(s,t))\mapsto (r,t)} \\
 \mathfrak{f} &= \int_R m_{(s,t) \mapsto t} d\lambda \\
 \mathfrak{g} &= m_{(r,s,t)\mapsto (s,t)} \\
 \mathfrak{h} &= a_{S,\int_R\omega_rd\lambda,(T^{'}_s),(\nu_s),(M^{'}_{(s,t)}) } \\
 \mathfrak{i} &= m_{((r,s),t) \mapsto (r,(s,t))} \\
 \mathfrak{j} &= a_{R\times S, \int_R \iota_r\omega_r d \lambda,(T^{''}_{(r,s)}),(\nu'_{r,s}),(M^{''}_{(r,s,t)})} \\
 \mathfrak{k} &= a_{R,\lambda,(S\times T)^{'},(\iota_r\int_S \iota_s\nu_sd\omega_r),(M^{'})} \\
 \mathfrak{l} &= \int_R a_{S ,\omega_r,(T^{'}),(\nu_s),(M^{'}_{(s,t)})}d\lambda \\
 \mathfrak{m} &= m_{(r,s)\mapsto s} \\
 \mathfrak{n} &= a_{R,\lambda,(S^{'}),(\omega_r),(\int_TM_td\nu_s)^{'}} 
\end{align*}
\end{minipage}
}
}
\label{diag7}
\end{figure}

\begin{figure}
\caption*{Diagram 8:}
\centering
\rotatebox{90}
{
\adjustbox{scale={0.47}{0.5}}
{
\begin{minipage}{52cm}

\begin{tikzcd}[ampersand replacement=\&]
	\& {\int_{\phi(R) \times \phi(S) \times \phi(T)}M^{''}_{(r,s,t)}d \int_R \iota_r (\int_S \iota_s \nu_s d\omega_r) d\lambda} \&\&\&\& {\int_{\phi(R) \times \phi(S) \times \phi(T)}M^{''}_{(r,s,t)}d \int_R \iota_r (\int_S \iota_s \nu_s d\omega_r) d\lambda} \&\&\&\&\&\& {\int_{\phi(R)} \int_{\phi(S) \times \phi(T)}M^{'} d \int_S \iota_s \nu_{s} d \omega_r d\lambda} \&\&\&\& \\
	\&\& 1 \&\&\&\&\& 2 \\
	{\int_{(R \times S) \times T}M^{''}_{(r,s,t)}d \int_R \iota_r (\int_S \iota_s \nu_s d\omega_r) d\lambda} \&\&\&\& {\int_{R \times (S \times T)}M^{''}_{(r,s,t)}d \int_R \iota_r (\int_S \iota_s \nu_s d\omega_r) d\lambda} \&\&\&\& {\int_R \int_{S \times T}M^{'} d \int_S \iota_s \nu_{s} d \omega_r d\lambda} \\
	\\
	\& {\int_{\phi(R) \times \phi(S)} \int_{\phi(T)}M_t  d\nu_{s}d \int_R \iota_r\omega_r d \lambda} \&\&\&\&\&\&\&\&\&\& {\int_{\phi(R)} \int_{\phi(S)} \int_{\phi(T)} M_t d\nu_s d \omega_r d\lambda} \&\&\&\& {} \\
	\&\&\&\& 5 \\
	{\int_{R \times S} \int_{T}M_t  d\nu_{s}d \int_R \iota_r\omega_r d \lambda} \&\&\&\&\&\&\&\& {\int_R \int_S \int_T M_t d\nu_s d \omega_r d\lambda}
	\arrow["{\mathfrak{a}}", from=1-2, to=1-6]
	\arrow[""{name=0, anchor=center, inner sep=0}, "{\mathfrak{b}}"{description, pos=0.8}, from=1-2, to=5-2]
	\arrow["{\mathfrak{c}}", from=1-6, to=1-12]
	\arrow[""{name=1, anchor=center, inner sep=0}, "{\mathfrak{d}}", from=1-12, to=5-12]
	\arrow["{\mathfrak{e}}"{description}, from=3-1, to=1-2]
	\arrow["{\mathfrak{f}}", from=3-1, to=3-5]
	\arrow[""{name=2, anchor=center, inner sep=0}, "{\mathfrak{g}}"{description}, shift left=3, from=3-1, to=7-1]
	\arrow["{\mathfrak{h}}"{description}, from=3-5, to=1-6]
	\arrow["{\mathfrak{i}}", from=3-5, to=3-9]
	\arrow["{\mathfrak{j}}"{description}, from=3-9, to=1-12]
	\arrow[""{name=3, anchor=center, inner sep=0}, "{\mathfrak{k}}"{description, pos=0.3}, from=3-9, to=7-9]
	\arrow["{\mathfrak{l}}", from=5-2, to=5-12]
	\arrow["{\mathfrak{m}}"{description, pos=0.6}, from=7-1, to=5-2]
	\arrow["{\mathfrak{n}}", from=7-1, to=7-9]
	\arrow["{\mathfrak{o}}"{description}, from=7-9, to=5-12]
	\arrow["3"', draw=none, from=2, to=0]
	\arrow["4"{description}, draw=none, from=3, to=1]
\end{tikzcd}

\vspace{2em}
\huge
\begin{align*} 
 \mathfrak{a} &= \la{Id} \\
\mathfrak{b} &= a^{'}_{\phi(R)\times \phi(S), \int_R \iota_r\omega_r d \lambda,(T^{''}_{(r,s)}),(\nu^{'}_{(r,s)}),(M^{''}_{(r,s,t)})} \\
\mathfrak{c} &= a^{'}_{\phi(R),\lambda,(\phi(S) \times \phi(T))^{'},(\iota_r\int_{\phi(S)} \iota_s\nu_sd\omega_r),(M^{'})} \\
\mathfrak{d} &= \int_{\phi(R)} a^{'}_{\phi(S) ,\omega_r,(T^{'}),(\nu_s),(M^{'}_{(s,t)})}d\lambda \\
\mathfrak{e} &= m_{((r,s),t) \mapsto (e_R(r),e_S(s),e_T(t))} \\
\mathfrak{f} &= m_{((r,s),t) \mapsto (r,(s,t))} \\
\mathfrak{g} &= a_{R\times S, \int_R \iota_r\omega_r d \lambda,(T^{''}_{(r,s)}),(\nu^{'}_{(r,s)}),(M^{''}_{(r,s,t)})} \\
\mathfrak{h} &= m_{(r,(s,t)) \mapsto (e_R(r),e_S(s),e_T(t))}\\
\mathfrak{i} &= a_{R,\lambda,(S\times T)^{'},(\iota_r\int_S \iota_s\nu_sd\omega_r),(M^{'})} \\
\mathfrak{j} &= (m \circ Tm)_{(r\mapsto e_R(r)),((s,t) \mapsto (e_S(s),e_T(t)))} \\
\mathfrak{k} &= \int_R a_{S ,\omega_r,(T^{'}),(\nu_s),(M^{'}_{(s,t)})}d\lambda \\
\mathfrak{l} &= a_{\phi(R),\lambda,(\phi(S)^{'}),(\omega_r),(\int_{\phi(T)} M_td\nu_s)^{'}} \\
\mathfrak{m} &= (m \circ Tm)_{(r,s)\mapsto (e_R(r),e_S(s)),(t \mapsto e_T(t))} \\
\mathfrak{n} &= a_{R,\lambda,(S^{'}),(\omega_r),(\int_TM_td\nu_s)^{'}} \\
\mathfrak{o} &= (m \circ Tm \circ T^2m)_{r \mapsto e_R(r),s \mapsto e_S(s),t \mapsto e_T(t)}=m_{r \mapsto e_R(r)} \circ \int_R (m\circ Tm)_{s \mapsto e_S(s),t \mapsto e_T(t)}d\lambda
\end{align*}
\end{minipage}
}
}
\label{diag8}
\end{figure}

\end{document}